\documentclass[11pt]{article}
\usepackage[margin=2cm]{geometry}
\usepackage[affil-it]{authblk}
\usepackage[utf8]{inputenc}
\usepackage{tikz}
\usepackage{keyval}
\usepackage{ifthen}
\usepackage{epsf}
\usepackage{psfrag}
\usepackage{wrapfig}
\usepackage{color,rotating,psfrag,amsfonts}
\usepackage{todonotes}
\usepackage[dvips]{epsfig}
\usepackage{adjustbox}

\usepackage{amsfonts}%
\usepackage{mathrsfs}%
\usepackage[title]{appendix}%
\usepackage{textcomp}%
\usepackage{manyfoot}%
\usepackage{booktabs}%

\usepackage{listings}%

\usepackage{graphicx} 
\usepackage{amsmath,amssymb,amsthm}
\usepackage{mathtools}
\usepackage{xcolor}
\usepackage{float}
\usepackage{diagbox}
\usepackage{subfigure}
\usepackage{multirow}
\usepackage{hyperref}
\usepackage{tikz-network}
\usepackage{scalefnt}
\usetikzlibrary{shapes.misc, positioning,  arrows.meta, shapes, shapes.multipart, positioning, shadows}

\usepackage{algorithm}%
\usepackage{algorithmicx}%
\usepackage{algpseudocode}%

\usepackage{varwidth}

\makeatletter
\newcounter{phase}[algorithm]
\newlength{\phaserulewidth}
\newcommand{\setphaserulewidth}{\setlength{\phaserulewidth}}
\newcommand{\phase}[1]{%
  \Statex\leavevmode\llap{\rule{\dimexpr\labelwidth+\labelsep}{\phaserulewidth}}\rule{\linewidth}{\phaserulewidth}
  
  \Statex\strut\refstepcounter{phase}{#1}
  \vspace{-1.25ex}\Statex\leavevmode\llap{\rule{\dimexpr\labelwidth+\labelsep}{\phaserulewidth}}\rule{\linewidth}{\phaserulewidth}}
\makeatother

\setphaserulewidth{.7pt}

\DeclarePairedDelimiter{\ceil}{\lceil}{\rceil}

\newtheorem{theorem}{Theorem}

\newtheorem{lemma}{Lemma}

\newtheorem{remark}{Remark}%

\newtheorem{definition}{Definition}%

\DeclareMathOperator*{\argmax}{arg\,max}
\DeclareMathOperator*{\argmin}{arg\,min}

\title{Cooperative games defined by multi-objective optimization in competition for  subsurface resources}

\date{}

\begin{document}

\author[1]{Per Pettersson\thanks{Corresponding author: per.pettersson@norceresearch.no}}

\author[2]{Sebastian Krumscheid}
\author[1,3]{Sarah Gasda}

\affil[1]{NORCE Norwegian Research Centre, N-5838 Bergen, Norway}
\affil[2]{Karlsruhe Institute of Technology (KIT), 76131 Karlsruhe, Germany}
\affil[3]{University of Bergen, N-5020 Bergen, Norway}

\maketitle

\begin{abstract}
We propose a novel decision making framework for forming potential collaboration among otherwise competing agents in subsurface systems. The agents can be, e.g., groundwater, CO$_2$, or hydrogen injectors and extractors with conflicting goals on a geophysically connected system. The operations of a given agent affect the other agents by induced pressure buildup that may jeopardize system integrity.
In this work, such a situation is modeled as a cooperative game 
where the set of agents is partitioned into disjoint coalitions that define the collaborations. The games are in partition function form with externalities, i.e., the value of a coalition depends on both the coalition itself and on the actions of external agents. 

We investigate the class of cooperative games where the coalition values are the total injection volumes as given by Pareto optimal solutions to multi-objective optimization problems subject to arbitrary physical constraints. For this class of games, we prove that the Pareto set of any coalition structure is a subset of any other coalition structure obtained by splitting coalitions of the first coalition structure. 
Furthermore, the hierarchical structure of the Pareto sets is used to reduce the computational cost in an algorithm to hierarchically compute the entire Pareto fronts of all possible coalition structures.

We demonstrate the framework on a pumping wells groundwater example, and nonlinear and realistic CO$_2$ injection cases, displaying a wide range of possible outcomes. Numerical cost reduction is demonstrated for the proposed algorithm with hierarchically computed Pareto fronts compared to independently solving the multi-objective optimization problems.
\end{abstract}

\section{Introduction}
\label{sec:intro}

Subsurface resources are assets with strategic and monetary value with a range of uses including energy extraction and energy storage, greenhouse gas storage, wastewater disposal, groundwater, and mining. The different users of subsurface resources often compete with each other over space and time, whether directly or indirectly by means of, e.g., induced pressure changes~\cite{Field_etal_18}. Due to increased populations and need for mitigation of climate change, the competition for subsurface resources (in particular CO$_2$ and hydrogen storage) and potential interaction between users is very likely to further increase in the future. Multi-agent systems models for subsurface resources, to be introduced in this paper, are relevant to account for the actions of independent agents with conflicting goals, operating on geophysically connected domains. These multi-agent models aim to ensure that all agents attain their respective goals in terms of value and risk avoidance, to the extent possible. 

The Intergovernmental Panel on Climate Change (IPCC) has  listed large-scale CO$_2$ storage as a necessary technology to reach the 1.5 $^{\circ}$C target by annual injection of 3-10 Gt during the coming  decades~\cite{IPCC_23}.
 Multi-site and basin-scale utilization of subsurface CO$_2$ storage resources are becoming feasible with maturing technology, while economic challenges remain and require further investigation~\cite{Krevor_etal_23}. In addition to the risk of directly jeopardizing reservoir integrity by pressure buildup~\cite{Rutqvist_etal_07}, CO$_2$ storage operations are likely to be affected by  natural gas storage and hydrocarbon extraction in hydraulically connected reservoirs. Safe and efficient utilization of large-scale CO$_2$ sites therefore calls for methods to optimize storage, where both physical and economic considerations are taken into account.

Very long time scales and large sites extending over hundreds of km, combined with complex multiscale physics, make numerical simulation of basin-scale CO$_2$ storage challenging. A remedy for the computational cost is offered by numerical models relying on physics-informed simplifications, e.g., vertical equilibrium models~\cite{Nilsen_etal_11}. Such models are essential in decision making, which requires both limited computational cost that allows repeated model evaluation, and sufficient accuracy that allows an informed decision.

Large-scale hydrogen storage in subsurface formations is predicted to play a crucial role in the energy transition and faces many of the same challenges as CO$_2$ storage. In addition, hydrogen will cyclically be injected in periods of relatively smaller energy demand, and extracted in periods of higher energy demand, leading to potentially critical stress changes due to pressure variation~\cite{Miocic_etal_23}. This can result in reactivation of faults, microseismic activity, and reservoir compaction, which may in turn lead to overburden rock subsidence, and formation of fractures acting as migration pathways~\cite{Tarkowski_UM_22}. Thus, the effects on agents operating on geologically connected domains can indeed be significant. While the challenges in large-scale hydrogen storage are related to those in large-scale CO$_2$ storage, they are even more complex, and call for multidisciplinary research for safe and efficient implementation~\cite{Heinemann_etal_21}.

Competition for subsurface resources due to physical interference effects also occurs in geothermal energy production, with multiple leases on a single unit. Unitization, where all lease owners coordinate their operations for optimized utilization, has repeatedly been proposed as a means to avoid negative effects, but has not been accepted as a worldwide standard~\cite{Goldstein_77, Tureyen_etal_15}.
Saline aquifers suitable for CO$_2$ storage may also have potential for geothermal energy extraction. Interaction between CO$_2$ storage and geothermal energy extraction on these sites need not be detrimental, as geothermal energy extraction may lower overpressurization~\cite{Randolph_etal_13}. Hence, there is potential for large gains by coordination of different subsurface activities on the same site.

In contrast to CO$_2$ and hydrogen storage, groundwater operations can successfully be described by single-phase models. While necessary to provide drinking water to nearby communities, these operations cause risks to agricultural yield, reduced river flow rates, and damaged wetland vegetation~\cite{Vink_Schot_02}. Hence, the benefits of accessible drinking water need to be weighed against the potential risks related to food production and ecological diversity.  The extensive literature on evolutionary algorithms for multi-objective optimization of groundwater resources reflect the challenges in competition for these resources: optimization of operating decisions in water resources systems~\cite{Takeuchi_Moreau_74}; multi-reservoir optimization with stochastic inflows~\cite{Wardlaw_Sharif_99}; well placement and rate optimization for groundwater pollution containment~\cite{Ritzel_etal_94}; cost and benefit optimization for drinking water~\cite{Vink_Schot_02}; comparison of evolutionary algorithms in groundwater multi-objective optimization~\cite{Hojjati_etal_18}, to give some examples.

The different settings described above, share the feature that there are often multiple independent agents (e.g., commercial companies), that may benefit from coordination of their actions by means of forming binding agreements, or at least benefit from considering what potential actions from external parties may affect themselves.
Cooperative game theory has received attention as a suitable framework to investigate problems related to coalition formation for rational decision making in multiagent systems.
A multiagent system consists of a finite number of agents that are free to form coalitions for cooperation to attain a value, assigned to the whole coalition and subsequently to be distributed among its members. The agents can be competitors and are typically only assumed to have an incentive to collaborate if they benefit from the collaboration.
Coalition formation can be divided into three stages~\cite{SANDHOLM_etal_99}: 1) Coalition structure generation, which amounts to
partitioning the agents into disjoint coalitions that coordinate their activities within but not between coalitions; 2) Optimization to achieve maximum (e.g., monetary) value within each coalition; and 3) Distribution of the attained value among the coalition members. 

The task of coalition structure generation is an NP-complete problem where the computational cost grows very fast with the number of agents and efficient algorithms are key, even under the common assumption that the value of every coalition is a-priori known~\cite{SANDHOLM_etal_99}. In the current setting of competition for subsurface resources, the number of agents will typically remain small, so the challenges related to the cost of coalition structure generation are not prohibitive. In contrast, determining (and distributing) the values of the coalitions is a challenging problem here, which can only be achieved to acceptable accuracy through computationally expensive numerical simulations of physical models. In addition, the values of the coalitions in a given partition of all agents may not be uniquely defined, even in situations where all agents have clearly defined goals. As an example, consider the additional simplifying assumption that all agents have comparable goals, e.g., all of them seek to maximize either CO$_2$ storage or hydrogen production subject to appropriate physical constraints. In this situation, injection or production cannot always be improved for one coalition without detriment to at least one other coalition. Multiple conflicting criteria typically lead to a range of outcomes that are all optimal, in the sense that none of them can be improved without negative impact on the others. In these cases, a decision maker, whether it is a reservoir engineer or a policy-maker with legislating power, should carefully make a decision about what outcomes to be considered the most relevant. This choice then defines the values of the associated cooperative game to be solved, i.e., given some additional criteria determining what is the desirable outcome of a game, identify what coalition structure should be formed based on the values computed.

Multi-criteria decision-making (MCDM) has been introduced to help systematically choosing among candidate solutions given technical, social, economic, and environmental criteria to be satisfied~\cite{Aruldoss_etal_13}. The criteria are often assumed to have been assigned weights in order of priority, and the goal of MCDM is to aggregate the frequently conflicting preferences. Informed decision making  typically requires input information from a complex numerical model that describe physical processes or other phenomena that impact the range of available options. The decision maker can be assumed to articulate preferences in different ways, depending on when and how the numerical model informs the decision process. A-priori articulation of preferences implies that the decision maker provides the relative order of preferences on the range of possible outcomes before evaluating the complex model. In a-posteriori preference articulation, the decision maker selects among the observed output after the model has been evaluated,
and in progressive preference articulation, the decision maker sequentially or continually provides information about preference to the model. 

MCDM methods can broadly be categorized into three classes that go under various names in the literature~\cite{Guitouni_Martel_98}: single synthesizing methods including multi-attribute utility theory methods~\cite{Keeney_Raiffa_76};  outranking methods~\cite{Roy_91}; interactive methods, including goal programming~\cite{Tamiz_etal_98} and Pareto front methods~\cite{Luce_Raiffa_57}. In the current work, we provide a general computational and exploratory framework for decision-making to outline a range of possibilities and will hence assume a-posteriori preference articulation. We also limit ourselves to directly comparable objectives among the agents, and believe that Pareto front methods are particularly well-suited for the settings considered with a small number of agents.

In this work, we form and analyze a class of cooperative games modeling competition about subsurface resources, where the games are "parameterized" by the non-unique solutions to multi-objective optimization problems involving numerical simulations of complex physical models. The parameterization implies that the games can have different qualitative properties, and they can even belong to different kinds of games, e.g., characteristic function games and partition function games, to be described in more detail below. We emphasize that the work presented here reverses the roles of game theory and optimization commonly considered in the literature. There is a long history of applying concepts from linear programming to solve games, as described in~\cite{Thie_Keough_11}. More recently, work has also focused on rewriting multi-objective optimization problems as different kinds of games where the players are the optimization objectives, c.f.~\cite{Meng_etal_10, Li_etal_21}.

The main contributions of this paper are as follows:
1) Propose a novel framework for decision-making where a cooperative game model has value functions defined by the solutions to multi-objective optimization problems; 2.1) Demonstrate the computational feasibility of the proposed methodology by presenting and 
proving
theoretical properties of the hierarchical characterization
of the Pareto sets to be computed that define the class of cooperative games; and 2.2) show numerical results for complex test cases where cost reduction is achieved by utilizing the hierarchical Pareto set properties; 3) Present exploratory results for decision-making, using the proposed framework. To achieve the third goal, we consider a-posteriori preference articulation to more broadly explore the range of different decision-making scenarios.

The paper is organized as follows. In Section~\ref{sec:methodology}, we introduce concepts from cooperative game theory, followed by a class of games defined by the solutions of multi-objective optimization problems. We prove theoretical hierarchical properties of the games where the objectives determining the value functions are linear in the values of the individual agents, but otherwise subject to complex physical constraints. The section is completed with an algorithmic workflow for numerically solving the corresponding cooperative games.
Numerical results for the proposed game-theoretical approach are presented in Section~\ref{sec:num_results}, where we first consider a linear groundwater problem which allows detailed illustration of game properties such as externalities, and stability properties of the coalition formation process. Then, two more complex test cases are investigated, both modeling CO$_2$ storage in the Bjarmeland formation in Barents Sea. We demonstrate that the theoretical properties of the games and multi-objective optimization problems can be used to increase the numerical efficiency of the proposed framework. Finally, we provide  a discussion to put the work into perspective and summarize the conclusions.

\section{Cooperative games with values from multi-objective optimization}
\label{sec:methodology}

\subsection{Cooperative Game Theory} 
In the following, we introduce terminology and notation from cooperative game theory to be used throughout the paper. This is similar to what can be found in textbooks on cooperative game theory~\cite{Chalkiadakis_etal_22}, and there are two main reasons for the choice of this abstract, mathematical style. First, the definitions are very general, and no adaptation to subsurface resources settings is needed or motivated. Second, while an abundance of more or less intricate results follow directly from the definitions given below, the definitions by themselves are quite natural and simple. We argue that everything needed in terms of understanding the game theory aspects of this paper is contained within these definitions. 

A cooperative game, to be  described in more detail below, consists of a number of  players, or \emph{agents}, that can choose to form binding agreements for collaboration and sharing value. A game also consists of a rule to assign values to all subgroups of agents. There are different solution concepts that identifies rational outcomes of the game, and sometimes also provides a systematic means to predict the outcome of the game. Let $A = \{a_1,a_2,\dots\}$ be the set of agents of size $|A|$. A \emph{coalition} $C$ is a non-empty subset of $A$. Let $\mathcal{C}_A := 2^{A}\setminus\emptyset$ denote the set of all possible coalitions among agents in $A$, where $2^{A}$ denotes the power set of $A$.
\begin{definition}[Coalition structure]
A \emph{coalition structure} $CS$ over $A$ is a partition of $A$ into disjoint coalitions, in the sense that $CS = \{C_1,C_2,\dots, C_k \}\subset\mathcal{C}_A$ for some $k \in \{1,2,\dots, |A| \}$ with $A = \cup_{j=1}^k C_j$ and $C_j\cap C_i = \emptyset$ for $i\not= j$. The set of all coalition structures over $A$ is denoted by $\Pi^{A}$.
\end{definition}
The number of coalition structures grows much faster than the number of coalitions, which itself increases as $2^{|A|}-1$, making exhaustive combinatorial explorations among all possible coalition structures infeasible for large problems. Depending on the character of the cooperative game describing the multiagent system, the coalition structure may or may not have an impact on the value of an individual coalition. We distinguish between two cases that will be relevant in what follows. 
\begin{definition}[Characteristic function game] A \emph{characteristic function game} (CFG) is given by a pair $(A,v)$, where $A$ is the set of agents, and the \emph{characteristic function} $v\colon\mathcal{C}_A\to \mathbb{R}$ maps each coalition $C$ to its value $v(C)$, which is independent of other (i.e., external) coalitions.
\end{definition}
It often matters what coalition structure a coalition belongs to. In such cases, it is useful to introduce the concept of embedded coalitions, i.e., the pair $(C,CS)$, where $C$ is a member of $CS$. The set of all embedded coalitions is denoted $\mathcal{EC}$. In the following, to simplify notation, we will often omit the reference to any particular $CS$, assuming that the fact that a coalition $C$ is embedded in $CS$ is clear from the context.
\begin{definition}[Partition function game]
A \emph{partition function game} (PFG) is given by a pair $(A,w)$, where $A$ is the set of agents, and the \emph{partition function} $w\colon\mathcal{EC}\to\mathbb{R}$ maps each coalition $C$ embedded in a coalition structure $CS \in \Pi^{A}$ to its value $w((C, CS))$, that also depends on coalitions in $CS$ that are different from $C$.
\end{definition}
In PFGs, the merging of two  coalitions can have a beneficial or adverse effect on the value of a given third coalition. The former case is referred to as positive externalities, and the latter, common in settings of finite resources, is referred to as negative externalities. More precisely, an \emph{externality} $\epsilon$ is defined as
\[
\epsilon(C,CS,CS') = w((C, CS)) - w((C, CS')),
\]
for any coalition $C$ in a coalition structure $CS = \{C, C_1 \cup C_2, CS'' \}$ and $CS' = \{ C, C_1, C_2, CS'' \}$, where $C_1 \neq C_2$ are some fixed coalitions, and $CS''$ is a coalition structure over $A \setminus ( C \cup C_1 \cup C_2) $. In the current work, where the value functions ($v$ or $w$) will be equal to any of the non-unique solutions of a multi-objective optimization problem, the more general concept of PFGs with externalities will be appropriate. Under certain conditions the problem simplifies to the special case of a CFG, which is by definition without externalities. We will use the notation $w$ and refer to it as a value function where it is assumed that it is either a partition function or a characteristic function depending on context.

A solution concept is a rule that determines rational outcomes of a game, i.e., selecting a coalition structure, and dividing the value among the coalition members (the third step of coalition formation). There are a few natural requirements for a feasible coalition structure. First, the coalition structure should be stable in the sense that no agent would benefit from leaving its coalition. 
Second, the coalition structure should also be fair so that each agent receives a payoff that corresponds to its contribution to the value of the coalition.
The \emph{payoff vector} $z \in \mathbb{R}^{|A|}$ assigns to each agent part of the value of the coalition they belong to, i.e., for each $(C, CS)$:
\[
\sum_{ \{ i | a_i \in C \} } z_i \leq w((C,CS)).
\]
The payoff is said to be \emph{efficient} if the inequality is replaced by equality in the above expression, i.e., if the entire value of a coalition is distributed to its members. While the value function and the payoff are closely related, the former assigns values to coalitions, and the latter (re)distributes values to the individual agents.

In a cooperative game with \emph{transferable utility}, the value of a coalition can be arbitrarily divided among the coalition members via the payoff vector, irrespective of their actual contribution to that value. In such games, the payoff of an agent may deviate from what would naturally be considered as its value. 
Assuming transferable utility adds to the complexity in the sense that payoff is negotiable and could be redistributed, leading to a wider range of possible outcomes of the game.

The efficiency of the \emph{grand coalition} (i.e., $C=A$) has often been assumed in previous literature and can be justified for CFGs, but not in general for PFGs~\cite{Hafalir_07}. That is, a coalition structure different from the grand coalition can generate larger values in the presence of externalities.

The concept of stability of a coalition structure is essential, but ambiguous in the context of PFGs. Hence, it deserves some further attention.
The set of all coalition structures of a CFG where no subset of agents has an incentive to deviate from the coalitions they belong to is called the \emph{core}~\cite{Gillies_59, Chalkiadakis_etal_22}. Sometimes in the literature, the core is defined as the set of payoffs and coalition structures where no coalition can benefit from leaving the grand coalition. The core may be empty. For PFGs, there is no unique definition of the core, since the value and payoff of a potential subcoalition depends on the actions taken by the members of external coalitions in the event of the given subcoalition breaking free from its allocation~\cite{Koczy_01}. In previous research, different general scenarios of coalition breakup have been considered, including the pessimistic approach where all remaining coalitions seek to penalize the payoff of the deviating agents~\cite{Lucas_67}, the optimistic approach where they seek to maximize the payoff of the deviators~\cite{Shapley_Shubik_66}, and the situation where an initial breakup always leads to all agents breaking free and forming singleton coalitions (i.e., no collaboration)~\cite{Chander_Tulkens_97}.

While the above examples illustrate the increased complexity and ambiguity in PFGs compared to CFGs, the current work will not suffer from this to any substantial extent for three reasons. First, there will be no default coalition structure from which agents may want to break free. Agents are assumed to start out independently and only form coalitions if they can benefit from that and not do any better in a different coalition. Hence, the concept of triggering reactions by deviating from a coalition does not apply. Second, the agents are assumed to only be interested in their own benefit and not actively try to work against or in favor of other agents unless it explicitly affects themselves. Third, we only analyze relatively small agent populations where all values and payoffs are already computed, so there is no need to make assumptions about unobserved behavior. Still, it is worth considering that the PFG models employed in the current work allow accounting for more complex dynamics, as indicated above. Despite the difficulties of introducing a unique definition of the core for general PFGs, the generic notion of the core as the set of stable coalition structures remain useful in its generality. Whether the stability by means of every agent's desire to remain within its coalition is guaranteed by legislation, signed contracts, superior payoff, or any other incentive, the core reduces the space of feasible coalition structures, among which eventually a single one should be selected that defines the outcome of the game.

In many cooperative games, it is of interest to investigate the total value of a coalition structure, known as \emph{social welfare}, and denoted $W(CS)$.
A coalition structure $CS^*$ that maximizes social welfare is defined as
\begin{equation}
\label{eq:soc_welfare}
CS^* = \argmax_{CS \in \Pi^{A}} W(CS), \quad \mbox{where} \quad W(CS) = \sum_{C \in CS} w((C,CS)),
\end{equation}
for a PFG, where it is assumed that the value of a coalition structure is the sum of the values of the coalitions. For a CFG, the social welfare-maximizing coalition structure is defined analogously, with $w((C, CS))$ replaced by $v(C)$.

Certain cooperative games have properties that make it superfluous to search the space of coalition structures for a winning solution. In superadditive games, the value of the union of any pair of coalitions is always greater than or equal to the sum of the values of each separate coalition. For these cases, the grand coalition of all agents will form, and the coalition structure generation problem becomes obsolete. We will encounter this situation as a special case in the setting of competition for finite subsurface resources, but emphasize that this will not be the case in general.

\subsection{Coalition values from multi-objective optimization}

The value functions in cooperative games are typically assumed directly available at unit cost in coalition structure generation problems~\cite{Rahwan_etal_15}. However, in the current work we consider the case where the value functions will not be available a-priori, but need to be computed as the solutions to multi-objective optimization problems instead, as is motivated in Sect.~\ref{sec:intro}. For a coalition $C$ in a coalition structure $CS$, the goal of multi-objective optimization  is to find a set of candidate decision variables $q \in \mathbb{R}^{N_{\text{dv}}}$, that simultaneously satisfy a number of possibly conflicting objectives $F_C^{CS}(q)$. The objective function is a parametrization of the value function in the sense that for every feasible $q$, the objective $F_C^{CS}(q)$ determines the corresponding value of $C \in CS$. However, the value function in the game theoretic setting is a function of the embedded coalition $(C, CS)$, hence the choice of separately defining an objective function of $q$. The value function for any embedded coalition $(C,CS)$ is then given by a corresponding objective for some carefully selected $q^*$:
\[
w((C,CS)) \equiv F_{C}^{CS}(q^{*}).
\]
For notational convenience and with some abuse of notation, we will henceforth drop the superscript $CS$, but it should be understood that the objectives pertain to some fixed coalition structure.

To determine an "optimal" decision variable vector $q$ for 
$|CS| \geq 2$ objectives, we seek the solution to the constrained multi-objective optimization problem (MOO)
\begin{equation}
\label{eq:moo}
\max_{q} \left(
F_{C_1}(q), F_{C_2}(q),\dots, F_{C_{|CS|}}(q)\right)\quad
\text{subject to}\quad g(q) \leq g_{\max},
\end{equation}
where $g\colon\mathbb{R}^{N_{\text{dv}}}\to\mathbb{R}^{N_{\text{con}}}$ is a nonlinear function of the decision variables, which encodes $N_{\text{con}}$ constraints via the condition $g(q) \leq g_{\max}$. We have excluded equality constraints, as they are not relevant for the applications we have in mind.
A solution to~\eqref{eq:moo} that is truly optimal for all objectives may not exist, and a more realistic goal is to instead find a balance of partially satisfying all objectives.
A \emph{Pareto efficient} solution has the property that no single objective can be improved without sacrificing at least one other objective. A candidate solution $q^{(1)}$ of a maximization problem is said to be (Pareto) \emph{dominated} by another candidate $q^{(2)}$ if 
\begin{equation}
\begin{aligned}
F_{C}(q^{(1)}) & \leq F_{C}(q^{(2)}) \mbox{ for all } C \in CS, \mbox{ and }\\
F_{C}(q^{(1)}) & < F_{C}(q^{(2)}) \mbox{ for at least one } C \in CS .
\end{aligned}
\end{equation}
In the case of minimization instead of maximization, the inequalities are reversed. 
The \emph{Pareto front} (PF) is the set of objective function vector values that correspond to the set of solutions (decision variables) that are not dominated by any other solutions. The Pareto (optimal) set (PS) is the set of  these non-dominated candidate solutions.

The MOO problem~\eqref{eq:moo} can be approximated directly by numerical MOO methods. Alternatively, it can first be rewritten as a set of single-objective optimization problems by means of the  \emph{weighted sum method} (WSM). In WSM, one seeks to maximize the weighted sum of the objective functions in~\eqref{eq:moo}, i.e.,
\begin{equation}
\label{eq:WSM}
F^{\text{WSM}}(q) = \sum_{C \in CS} \alpha_{C} F_{C}(q),   
\end{equation}
with weights $\alpha_C\ge 0$ for all $C\in CS$ and $\sum_{C \in CS} \alpha_{C}=1$. If the weights are strictly positive, then the WSM is sufficient for Pareto optimality~\cite{Zadeh_63}, a fact that will later be used in both theoretical and numerical results and is summarized in the following lemma.
\begin{lemma}
\label{lemma:WSM_PF}
The unique solution to the weighted sum method with arbitrary but fixed strictly positive weights applied to the objective functions of an MOO is always a Pareto optimal solution to the MOO.
\end{lemma}
The maximization of the WSM formulation is furthermore a necessary condition in the case of convex problems, i.e., the constraints encoded in $g$ and objective functions $F_C$ are all convex~\cite{Geoffrion_68}.
In particular, since WSM captures the convex parts of the PF, if the true PF is indeed convex, then the WSM method can capture the full PF.

\subsection{Pareto front point selection}
\label{sec:pt_select}
Once the PF has been identified, it is desirable to use systematic means to select the preferred solution, or a reduced subset of candidate solutions that can be subsequently investigated, so that a decision maker can determine a final unique selection. Methods to perform PF selection include the technique for order of preference by similarity to ideal solution (TOPSIS) where the solution is sought to have the smallest possible distance to a positive ideal solution and largest possible distance to a negative ideal solution~\cite{Tzeng_Huang_81}; the class of ELECTRE~\cite{Roy_68, Roy_91} and PROMETHEE~\cite{Brans_etal_86} methods based on outranking approaches; the widely used analytic hierarchy process~\cite{Saaty_90} that can be applied to PFs~\cite{Ayadi_etal_17}. A review of quantitative methods for PF point selection, including detailed mathematical expressions, is given in~\cite{Wang_Rangaiah_17}.
It should come as no surprise that there is some overlap in these methods with those referenced in the description of MCDM methods in the Introduction. We emphasize however that the  focus in this subsection is not decision-making with respect to the global problem of choosing collaboration partners by determining a winning coalition structure, but on finding ways to motivate or single out a certain outcome among viable alternatives for any given coalition structure. 
A popular class of methods is utopia 
 point methods~\cite{Yu_13}: select the point at the PF that in some metric minimizes the distance to the Utopia point $F^*\in\mathbb{R}^{|CS|}$, whose components $F_{C}^*$ are defined by
\[
F_{C}^* = \max_{q} \{ F_{C}(q) | g(q) \leq g_{\max} \}, \quad \forall C \in CS.
\]
A solution that coincides with the utopia point does not exist in general, and instead one may seek a Pareto solution $q$ that minimizes the distance between  $F(q)$ and the utopia point in a weighted $L_p$ norm with $1\leq p<\infty$~\cite{Lu_etal_11}, i.e., 
\begin{equation}
\label{eq:utopia_norm}
F_C^{\beta} \equiv F_{C}(q_{\beta}^\ast) \quad\text{where}\quad q^{*}_{\beta} = \argmin_{q \in \text{PS}} \left( \sum_{C \in CS} |\beta_C (F_{C}(q) -F_{C}^{*})|^{p}\right)^{1/p},
\end{equation}
where $\{\beta_{C}\}_{C \in CS}$ are positive weights chosen by the decision maker. Setting $p=1$ yields solutions on the convex hull of the PF. Increasing $p$ widens the set of points of the PFs that can be chosen by varying the weights, where $p=\infty$ includes the full PF~\cite{Kasprzak_Lewis_01}.

In this work, we consider a set of PFs, each representing the non-dominated solutions to a game with coalition structure $CS \in \Pi^{A}$ and with $|CS|$ objectives. The $|\Pi^{A}|$ MOOs are distinct but related (indeed, they model the same constrained physical setup but with different coalition structures), and thus the weights $\{\beta_{C}\}_{C \in CS}$ should be chosen consistently between the MOOs. We assume that each agent $a\in A$ receives a weight, $\tilde{\beta}_{a}$ that does not change with the coalition structure. The coalition weights are then the (unweighted) sums of individual agent weights, i.e., for coalition structure $CS$ we set 
\begin{equation}
\label{eq:weighted_agents}
\beta_{C} = \sum_{ 
\{ a | a  \in C \}
} \tilde{\beta}_{a}, \quad \forall C \in CS. 
\end{equation}

As an alternative PF selection criterion, we propose to choose the non-dominated solution that provides the best payoff $z_a$ for some agent $a$. Recall that throughout this work, we assume that the payoff of an agent is equal to its value, i.e., $z_a = w((\{ a\},CS))$, where the value is determined by the objective function evaluated at the selected point from the PS. To highlight the dependence on the particular point $q$ of the PS, we will use the notation $z_a \equiv z_a(q)$.
Mathematically, we then choose $F_C^{\text{max }a}$ defined by 
\begin{equation}
\label{eq:favor_ind_agent}
 F_C^{\text{max }a} \equiv F_{C}(q_a^\ast)\quad\text{where}\quad z_{a} = \max_{q \in \text{PS}} z_{a}(q)\quad\text{and}\quad q_a^\ast = \argmax_{q \in \text{PS}} z_{a}(q)
\end{equation}
Note that this criterion is not equivalent to setting $\tilde{\beta}_{a}$ to some large number in~\eqref{eq:weighted_agents}. A large value of $\tilde{\beta}_{a}$ favors the entire coalition that $a$ belongs to, while the formulation~\eqref{eq:favor_ind_agent} only favors the payoff for the agent $a$ itself.

\subsection{Games defined by objectives linear in agents and complex constraints}

Every combination of selections of a PF point from each MOO defines a cooperative game, and all possible combinations define a class of games. 
For the applications we are considering here, the constraint function can be highly nonlinear, and its evaluation requires a numerical simulation of a physical problem. In contrast, the objective function is often a linear or almost linear function of the value functions of the individual agents. It is a relatively weak assumption that the value of a coalition is the sum of the values of its subcoalitions.
Hence, a general and highly relevant class of games are defined by the MOO~\eqref{eq:moo} with constraints that can be arbitrarily complex, as given by, e.g., a complex physical model, as long as they are the same for all coalition structures. That is, the feasible solution space can be non-convex but assumed to remain the same for all $CS \in \Pi^{A}$. The objectives are assumed to be weighted sums of the values of the individual agents in the coalitions, which can be arbitrary different functions of the design variables. Specifically, the objectives in~\eqref{eq:moo} are given by
 \begin{equation}
 \label{eq:gen_obj_lin_agents}
 F_{C}(q)= \sum_{a \in C} \gamma_{a}F_{a}(q), \quad \forall C \in CS .
 \end{equation}

 Next, we present some theoretical results that hold for all MOO problems with objectives given by~\eqref{eq:gen_obj_lin_agents}. First, note that, if $q^*$ is feasible for one $CS$, it is feasible for all $CS$ since the constraint function is the same for all $CS$, no matter its complexity. Unless one considers transferable utility, it is natural to equate the payoff of each agent with the corresponding values directly determined by $q^{*}$. 

\begin{theorem}   \label{thm:hierarch_PS}
Let $A$ be the set of agents and denote by $\Gamma$ be the class of cooperative games for $A$ with coalition values being Pareto optimal solutions to the MOO~\eqref{eq:moo}, where the objectives are given by~\eqref{eq:gen_obj_lin_agents}. 
Furthermore, let $CS \neq CS'$ be any two coalition structures such that $CS$ is a refinement of $CS'$, in the sense that there exists a subset $S\subseteq CS$ so that $C' = \cup_{C\in S} C$ for at least one coalition $C' \in CS'$. Then it holds for the Pareto sets corresponding to those coalition structures that $PS_{CS'} \subseteq PS_{CS}$.
\end{theorem}

\begin{proof}
Let the coalition structure $CS$ be a refinement of the coalition structure $CS'$. 
 By the linearity of the objectives in the agents' values as defined in 
 Eq.~\eqref{eq:gen_obj_lin_agents} and the fact that coalitions in $CS$ are disjoint, we have that for $C'\in CS'$
 \[
 F_{C'}(q) = \sum_{a \in C'}\gamma_{a}F_a(q) 
 = \sum_{C\in S} \sum_{a \in C} \gamma_{a}F_a(q)
 = \sum_{C\in S} F_{C}(q),
 \]
showing that the objective of a larger coalition is always the sum of the objectives of its subcoalitions.

To show $PS_{CS'} \subseteq PS_{CS}$, let $q'$ be a member of $ PS_{CS'}$ so that it is any Pareto optimal solution corresponding to $CS'$ and we need to show that $q'$ is also a member of $PS_{CS}$. We will prove this by assuming the opposite, which will lead to a contradiction. For this, assume that $q'$ is not a Pareto optimal solution corresponding to $CS$. Hence, $q'$ is dominated by some solution $q^*$ with respect to $CS$ so that $F_{C}(q^*) \geq F_{C}(q')$ for all $C \in CS$ (and $F_{C}(q^*) > F_{C}(q')$ for at least one $C \in CS$). However, from the Pareto optimality of $q'$ with respect to $CS'$, for some $C' \in CS'$ we have $F_{C'}(q') > F_{C'}(q^*)$, which can be written as 
\[
F_{C'}(q') = \sum_{C \in S} F_{C}(q') > \sum_{C\in S} F_{C}(q^*) = F_{C'}(q^*).
\]
This leads to a contradiction of the dominance of $q'$ over $q^*$, hence there cannot exist a $q^*$ that dominates $q'$ with respect to $CS$.
Thus, if $q'$ is a member of the Pareto set $PS_{CS'}$, then it is a member of the Pareto set $PS_{CS}$, as claimed.
 \end{proof}

\begin{remark}
In general, the Pareto sets are proper subsets in the above Theorem~\ref{thm:hierarch_PS}, i.e., $ PS_{CS'} \subset PS_{CS}$. To see this, we observe that by Lemma~\ref{lemma:WSM_PF}, one can apply the WSM to both $CS'$ and $CS$ so that the unique optima of, respectively,
 \begin{align}
  F^{\textup{WSM}}_{CS'}(q) &\equiv  \sum_{C'  \in CS'}  \alpha_{C'} F_{C'}(q) = \sum_{C'  \in CS'}  \alpha_{C'}  \sum_{C \in S}  F_{C}(q),
\\
    F_{CS}^{\textup{WSM}}(q) &\equiv \sum_{C  \in CS}  \alpha_{C} F_{C}(q) = \sum_{C' \in CS'}  \sum_{C\in S} \alpha_{C} F_{C}(q),
 \end{align}
 belong to $PF_{CS'}$ and $PF_{CS}$, respectively. The expressions are identical (and hence guarantee the same optimal $q^*$) only if $\alpha_{C}= \alpha_{C'}$ for all $C\in S$  
 and for all $C' \in CS'$.
 By choosing $\alpha_{C_i} \neq \alpha_{C_j}$ for some $C_i, C_j \subset C'$, the WSM yields Pareto optimal solutions for $CS$ that are not in general Pareto optimal solutions of $CS'$, and $PS_{CS'}\subset PS_{CS}$. Note that this is not a proof. For instance, in the exceptional case where the feasible solution set due to very limiting constraints contains a single point only, so that all Pareto sets consists of the same single point, all WSM weights yield the same solution, and $PS_{CS} = PS_{CS'}$.
\end{remark}
A few consequences of Theorem~\ref{thm:hierarch_PS} are noteworthy. 
\begin{remark}
The unique solution that belongs to all Pareto sets of the theorem always maximizes the social welfare~\eqref{eq:soc_welfare}, and it is the solution closest to the utopia point in the unweighted norm. Provided that maximization of social welfare (e.g., optimal utilization of subsurface resources) is the criterion used to determine the outcome of the game, there will be no incentive for agents to form any kind of collaboration. 
\end{remark}
\begin{remark}
To obtain the social-welfare maximizing solution, it is sufficient to  solve only a single constrained single-objective optimization problem to solve the coalition structure generation problem, as opposed to a potentially very large numbers of constrained MOO problems. Depending on the problem of interest, this may imply orders of magnitude numerical cost reduction.
\end{remark}
\begin{remark}
For any number of agents, it is in principle sufficient to compute a single PF corresponding to the singleton coalition structure and then obtain all the PFs corresponding to the remaining coalition structures via postprocessing of the PS of the first PF.
\end{remark}

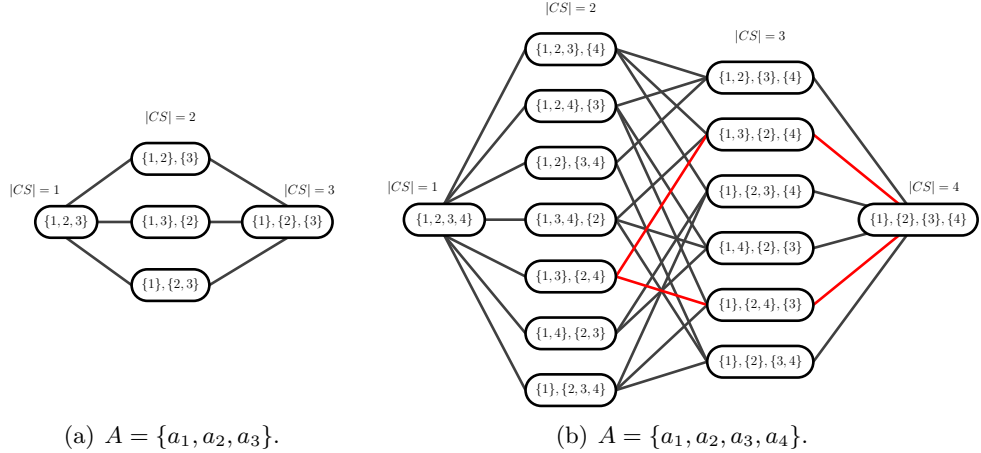
\begin{figure}[h]
\centering
\subfigure[$A=\{ a_1, a_2, a_3 \}$.]  
{ 
\begin{tikzpicture}[thick,scale=0.42, every node/.style={transform shape}]
\Text[fontsize=\large, x=6.7, y=0]{$|CS|=3$}
\Vertex[shape=rounded rectangle, color=white, size=1.0, fontsize=\large, x=6.0, y=-1, label={ \ $\{1\}, \{ 2 \}, \{ 3\} $ \ }, style={minimum width=3.0cm}]{one-two-three}
\Text[fontsize=\large, x=2.3, y=2.3]{$|CS|=2$}
\Vertex[shape=rounded rectangle, color=white, size=1.0, fontsize=\large, x=2.3, y=1, label={ \ $\{1,  2 \}, \{ 3\} $ \ }, style={minimum width=2.6cm}]{onetwo-three}
\Vertex[shape=rounded rectangle, color=white, size=1.0, fontsize=\large, x=2.3, y=-1, label={ \ $\{1,  3 \}, \{ 2\} $ \ }, style={minimum width=2.6cm}]{onethree-two}
\Vertex[shape=rounded rectangle, color=white, size=1.0, fontsize=\large, x=2.3, y=-3, label={ \ $\{1\}, \{ 2,  3\} $ \ }, style={minimum width=2.6cm}]{one-twothree}
 \Text[fontsize=\large, x=-2, y=0]{$|CS|=1$} 
 \Vertex[shape=rounded rectangle, color=white, size=1.0, fontsize=\large, x=-1, y=-1, label={ \ $\{1,  2,  3\} $ \ }, style={minimum width=2.1cm}]{onetwothree}

 \Edge[lw=1, label={}, fontsize=\large](one-two-three.north)(onetwo-three.east)
 \Edge[lw=1, label={}, fontsize=\large](one-two-three)(onethree-two)
 \Edge[lw=1, label={}, fontsize=\large](one-two-three.south)(one-twothree.east)
    
 \Edge[lw=1, label={}, fontsize=\large](onetwothree.north)(onetwo-three.west)
 \Edge[lw=1, label={}, fontsize=\large](onetwothree)(onethree-two)
 \Edge[lw=1, label={}, fontsize=\large](onetwothree.south)(one-twothree.west)

  \Text[fontsize=\large, x=-2, y=-6.8]{}  
    
\end{tikzpicture}
}
~
\subfigure[$A=\{ a_1, a_2, a_3, a_4\}$.]  
{ 
\begin{tikzpicture}[thick,scale=0.42, every node/.style={transform shape}]
\Text[fontsize=\large, x=16.5, y=1.0]{$|CS|=4$}
\Vertex[shape=rounded rectangle, color=white, size=1.0, fontsize=\large, x=16, y=0, label={ \ $\{1\}, \{ 2 \}, \{ 3\}, \{ 4 \} $ \ }, style={minimum width=3.9cm}]{one-two-three-four}

\Text[fontsize=\large, x=11, y=5.8]{$|CS|=3$}
\Vertex[shape=rounded rectangle, color=white, size=1.0, fontsize=\large, x=11, y=4.5, label={ \ $\{1,  2 \}, \{ 3\}, \{ 4\} $ \ }, style={minimum width=3.5cm}]{onetwo-three-four}
\Vertex[shape=rounded rectangle, color=white, size=1.0, fontsize=\large, x=11, y=2.7, label={ \ $\{1,  3 \}, \{ 2\}, \{ 4\} $ \ }, style={minimum width=3.5cm}]{onethree-two-four}
\Vertex[shape=rounded rectangle, color=white, size=1.0, fontsize=\large, x=11, y=0.9, label={ \ $\{1\}, \{ 2,  3\}, \{ 4\} $ \ }, style={minimum width=3.5cm}]{one-twothree-four}

\Vertex[shape=rounded rectangle, color=white, size=1.0, fontsize=\large, x=11, y=-0.9, label={ \ $\{1, 4\}, \{ 2 \}, \{ 3\} $ \ }, style={minimum width=3.5cm}]{onefour-two-three}
\Vertex[shape=rounded rectangle, color=white, size=1.0, fontsize=\large, x=11, y=-2.7, label={ \ $\{1\}, \{ 2, 4 \}, \{ 3\} $ \ }, style={minimum width=3.5cm}]{one-twofour-three}
\Vertex[shape=rounded rectangle, color=white, size=1.0, fontsize=\large, x=11, y=-4.5, label={ \ $\{1\}, \{2\}, \{ 3,  4\} $ \ }, style={minimum width=3.5cm}]{one-two-threefour}

 \Text[fontsize=\large, x=5, y=6.7]{$|CS|=2$}    
 \Vertex[shape=rounded rectangle, color=white, size=1.0, fontsize=\large, x=5, y=5.4, label={ \ $\{1,  2,  3\}, \{ 4\} $ \ }, style={minimum width=3.0cm}]{onetwothree-four}
 \Vertex[shape=rounded rectangle, color=white, size=1.0, fontsize=\large, x=5, y=3.6, label={ \ $\{1,  2,  4\}, \{ 3\} $ \ }, style={minimum width=3.0cm}]{onetwofour-three} 
 \Vertex[shape=rounded rectangle, color=white, size=1.0, fontsize=\large, x=5, y=1.8, label={ \ $\{1,  2\}, \{ 3,  4\} $ \ }, style={minimum width=3.0cm}]{onetwo-threefour} 
  \Vertex[shape=rounded rectangle, color=white, size=1.0, fontsize=\large, x=5, y=0, label={ \ $\{1,  3, 4\}, \{ 2\} $ \ }, style={minimum width=3.0cm}]{onethreefour-two} 
  \Vertex[shape=rounded rectangle, color=white, size=1.0, fontsize=\large, x=5, y=-1.8, label={ \ $\{1,  3\}, \{ 2, 4\} $ \ }, style={minimum width=3.0cm}]{onethree-twofour} 
    \Vertex[shape=rounded rectangle, color=white, size=1.0, fontsize=\large, x=5, y=-3.6, label={ \ $\{1,  4\}, \{ 2, 3 \} $ \ }, style={minimum width=3.0cm}]{onefour-twothree} 
 \Vertex[shape=rounded rectangle, color=white, size=1.0, fontsize=\large, x=5, y=-5.4, label={ \ $\{1\}, \{ 2, 3,  4 \} $ \ }, style={minimum width=3.0cm}]{one-twothreefour}
 
  \Text[fontsize=\large, x=0, y=1.0]{$|CS|=1$}    
 \Vertex[shape=rounded rectangle, color=white, size=1.0, fontsize=\large, x=1, y=0, label={ \ $\{1,  2,  3, 4\} $ \ }, style={minimum width=2.7cm}]{onetwothreefour}

 \Edge[lw=1, label={}, fontsize=\large](one-two-three-four)(onetwo-three-four.east)
 \Edge[lw=1, label={}, color=red, fontsize=\large](one-two-three-four)(onethree-two-four.east)
 \Edge[lw=1, label={}, fontsize=\large](one-two-three-four)(one-twothree-four.east)
  \Edge[lw=1, label={}, fontsize=\large](one-two-three-four)(onefour-two-three.east)
 \Edge[lw=1, label={}, color=red, fontsize=\large](one-two-three-four)(one-twofour-three.east)
 \Edge[lw=1, label={}, fontsize=\large](one-two-three-four)(one-two-threefour.east)
    

 \Edge[lw=1, label={}, fontsize=\large](onetwo-three-four.west)(onetwothree-four.east)
 \Edge[lw=1, label={}, fontsize=\large](onetwo-three-four.west)(onetwo-threefour.east)
  \Edge[lw=1, label={}, fontsize=\large](onetwo-three-four.west)(onetwofour-three.east)

 \Edge[lw=1, label={}, fontsize=\large](onethree-two-four.west)(onetwothree-four.east)
 \Edge[lw=1, label={}, fontsize=\large](onethree-two-four.west)(onethreefour-two.east)

\Edge[lw=1, label={}, fontsize=\large](one-twothree-four.west)(one-twothreefour.east)
\Edge[lw=1, label={}, fontsize=\large](one-twothree-four.west)(onefour-twothree.east)
\Edge[lw=1, label={}, fontsize=\large](one-twothree-four.west)(onetwothree-four.east)

\Edge[lw=1, label={}, fontsize=\large](onefour-two-three.west)(onetwofour-three.east)
\Edge[lw=1, label={}, fontsize=\large](onefour-two-three.west)(onefour-twothree.east)
\Edge[lw=1, label={}, fontsize=\large](onefour-two-three.west)(onethreefour-two.east)

 \Edge[lw=1, label={}, fontsize=\large](one-twofour-three.west)(one-twothreefour.east)

  \Edge[lw=1, label={}, fontsize=\large](one-twofour-three.west)(onetwofour-three.east)
  
  \Edge[lw=1, label={}, fontsize=\large](one-two-threefour.west)(one-twothreefour.east) 
   \Edge[lw=1, label={}, fontsize=\large](one-two-threefour.west)(onetwo-threefour.east)
   \Edge[lw=1, label={}, fontsize=\large](one-two-threefour.west)(onethreefour-two.east)  

  \Edge[lw=1, label={}, color=red, fontsize=\large](onethree-two-four.west)(onethree-twofour.east)
  \Edge[lw=1, label={}, color=red, fontsize=\large](one-twofour-three.west)(onethree-twofour.east)  
 
  \Edge[lw=1, label={}, fontsize=\large](onetwothreefour.north)(onetwothree-four.west)
 \Edge[lw=1, label={}, fontsize=\large](onetwothreefour.north)(onetwofour-three.west)
 \Edge[lw=1, label={}, fontsize=\large](onetwothreefour.north)(onetwo-threefour.west)
  \Edge[lw=1, label={}, fontsize=\large](onetwothreefour.east)(onethreefour-two.west)
 \Edge[lw=1, label={}, fontsize=\large](onetwothreefour.south)(onethree-twofour.west)
 \Edge[lw=1, label={}, fontsize=\large](onetwothreefour.south)(onefour-twothree.west)
  \Edge[lw=1, label={}, fontsize=\large](onetwothreefour.south)(one-twothreefour.west)
    
\end{tikzpicture}
}
\caption{Coalition structure graphs for (a) three agents; and (b) four agents, respectively. Following the edges colored in red from a target coalition structure $\{ \{a_1, a_3\}, \{a_2, a_4\} \}$ shows the coalition structures whose corresponding Pareto sets are supersets of the target Pareto set, i.e., $\{\{a_1, a_3\}, \{a_2\}, \{a_4\}\}$, $\{\{a_1\}, \{a_2, a_4\}, \{a_3\} \}$, and $\{\{a_1\}, \{a_2\}, \{a_3\}, \{a_4\}\}$.}
\label{fig:coal-struct-graphs}
\end{figure}

As an illustration and to appreciate the meaning of Theorem~\ref{thm:hierarch_PS}, first consider the coalition structure graph for a set of agents, obtained as follows. Every node in the graph represents a coalition structure, and the nodes are organized in levels according to the number of coalitions they contain. The first level contains the grand coalition only, the second level all coalition structures with two coalitions, and so on. There is an undirected edge between two nodes only if the coalition structure on the coarser level (i.e., the one with fewer number of coalitions) can be obtained by splitting exactly one coalition into two coalitions to obtain the coalition structure on the finer level.
The corresponding coalition structure graphs for three and four agents are shown in Figure~\ref{fig:coal-struct-graphs}. In view of a coalition structure graph, the interpretation in terms of Theorem~\ref{thm:hierarch_PS} is that the PS of any given coalition structure is a subset of the PS of the coalition structures corresponding to nodes that are connected by following the edges to the right in the graph. 
For instance, if we want to obtain the PS corresponding to the coalition structure $\{ \{a_1, a_3\}, \{a_2, a_4\} \}$ it is a subset of any of the coalition structures PS that can be reached by the edges colored in red in Figure~\ref{fig:coal-struct-graphs}, i.e., $\{\{a_1, a_3\}, \{a_2\}, \{a_4\}\}$, $\{\{a_1\}, \{a_2, a_4\}, \{a_3\} \}$, and $\{\{a_1\}, \{a_2\}, \{a_3\}, \{a_4\}\}$.
While the situation is relatively simple for the case of three agents, the graph with four agents indicates that the hierarchical structure of the PS can be quite intricate for problems with larger number of agents. The coalition structure graph can be used both to estimate the numerical cost reduction by employing the hierarchical structure in the development of a hierarchically coupled MOO algorithm, and as a bookkeeping structure in its implementation.


\subsection{A workflow for decision-making in cooperative games defined by MOO}
\label{sec:algorithms}

\begin{figure}[h]
\centering
    \begin{tikzpicture}[scale=0.9, transform shape,every text node part/.style={align=center},
node distance = 5mm and 7mm,
   arr/.style = {-Triangle,  thick},
   box/.style = {rectangle, draw, semithick,
                 minimum height=9mm, minimum width=17mm,
                 fill=white} 
                        ]      
                        
  \tikzstyle{decision} = [diamond, minimum width=5cm, minimum height=1cm, text centered, draw=black, semithick, aspect=2, inner xsep=0mm]
   
\node (n1) [box, rounded corners] {\textbf{Initialize Class of Cooperative Games} \\
\small{
\begin{varwidth}{\linewidth}
            \begin{algorithmic}
                \State $A \leftarrow$ Set of agents
                \State $\Pi^A \leftarrow$ Set of all coalition structures
                \State $w(C, CS; q) \forall C \in CS \in \Pi^{A}$ Value function
                \State $q \leftarrow$ Design variables
                \State $g(q) \leftarrow$ Physical and economic constraints
        \end{algorithmic}%
        \end{varwidth}
        }
        };
       
\node (n3) [decision, below=of n1] {Theorem~\ref{thm:hierarch_PS} applies?};

\node (n4) [box, below left=of n3] {\textbf{Hierarchical MOO} \\ 
\small{
\begin{varwidth}{\linewidth}
  Algorithm~\ref{alg:hier_MOO}
        \end{varwidth}
        }
    
};
\node (n5) [box, below right=of n3] {\textbf{Non-nested MOO} \\ 
\small{
\begin{varwidth}{\linewidth}
Algorithm~\ref{alg:standard_MOO}
        \end{varwidth}
        }
};

\node (n6) [box, below right=of n4] {\textbf{Pareto fronts}};

\node (n7) [draw, trapezium,trapezium left angle=70,trapezium right angle=-70, semithick, below=of n6] {\textbf{Decision making 1}\\ Apply Pareto front \\ selection criterion};
\node (n8) [draw, trapezium,trapezium left angle=70,trapezium right angle=-70, semithick, below=of n7] {\textbf{Decision making 2} \\ Apply solution concept to game};

\draw[arr]   (n1) -- (n3);

\draw[arr]   (n3) -| node [near start, above] {Yes} (n4.north);
\draw[arr]   (n3) -| node [near start, above] {No} (n5.north);

\draw[arr]   (n4.south) |-  (n6.west);
\draw[arr]   (n5.south) |- (n6.east);

\draw[arr]   (n6) -- (n7);
\draw[arr]   (n7) -- (n8);                          

\end{tikzpicture}

\caption{Overview of full workflow with class of cooperative games with value functions defined by the solutions to multi-objective optimization problems with physical constraints.}
\label{fig:flowchart_algo}
\end{figure}
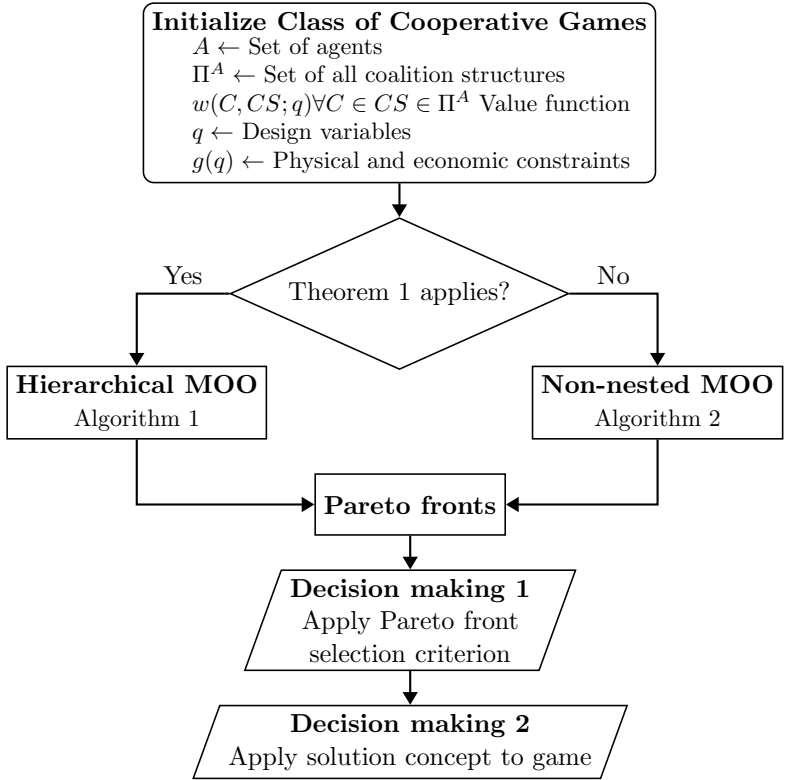

We are now ready to present a full workflow for decision-making using the presented concepts from cooperative game theory and MOO, as illustrated in Figure~\ref{fig:flowchart_algo}. First, the cooperative game needs to be defined in terms of agents, objectives that will determine the value function, design variables controlled by the agents or external parties, and a physical model constraining actions of the agents and the values that can be attained. This phase defines the problem, i.e., the class of cooperative games, but it is not computationally intensive since the model has not yet been numerically evaluated. 

Next, the computationally expensive phase of computing the value functions by MOO follows. If Theorem~\ref{thm:hierarch_PS} applies, then the computational cost can be significantly reduced, as will be described in more detail later in this section. The PFs that result from the MOO constitute all possible coalition values, for all coalition structures of the class of games we investigate. With the a-posteriori preference articulation considered in this work, a decision maker can subsequently apply different PF point selection criteria to thus obtain individual members of the class of games. We emphasize that the criteria are not limited to the examples proposed in Section~\ref{sec:pt_select}. To reach a final decision about subsurface resource utilization, the decision maker needs to select and apply a solution concept. The two decision-making steps, PF selection and choice of solution concept, can be varied to explore the solution space more broadly without the need for additional model evaluations of the computationally expensive components of the framework.

Next, we discuss in more detail the numerical solution of the set of MOOs,
which is the computationally most intensive component of the proposed workflow. To be clear, we investigate how the hierarchical structure of the PS can be utilized to more efficiently obtain all Pareto solutions, rather than focusing on the choice of numerical method for any single MOO problem, for which there is already a rich literature, c.f.~\cite{Cui_etal_17}.
By virtue of Theorem~\ref{thm:hierarch_PS}, the solution to the grand coalition optimization problem is a member of the PSs of all the other coalition structures, and more generally, the PS of any coalition structure $CS$ is a  subset of the PS of any other refinement of $CS$. It follows that it is sufficient to solve only the multi-objective optimization problem corresponding to the singleton coalition structure to desired accuracy, and the solutions to all other coalition structures can be obtained by post-processing at very modest computational cost. Since we start from the most refined coalition, which is often placed at the top level of the coalition structure graph, and then obtain the solutions for all lower levels, we will refer to this as the \emph{top-down} approach. Depending on the choice of multi-objective optimization method, a subset or approximation $q'$ of the set of non-dominated solutions $q^{*}$ can be used to more efficiently find better approximations of the full set of non-dominated solutions. We use the following notation for the MOO with available (partial) approximations $q'$ of the solution:
\[
q^{*}=\argmax_{q} (F_{C_1}(q),\dots, F_{C_{|CS|}}(q) | q' ).
\]
Hence, rather than starting directly with computing the singleton coalition structure PF, one can start with the simpler problem of computing, e.g., the grand coalition solution, and use that as an initial guess for the coalition structures with exactly two coalitions. Then, these solutions, that may be a subset of all non-dominated solutions for that coalition structure, can be used as an initial guess for the next level of refinement.  These considerations can be combined into an algorithm with three stages, as presented in Algorithm~\ref{alg:hier_MOO}. The first stage is optional, and hierarchically computes approximations to some Pareto optimal solutions. Starting with the coarsest coalition structure (the grand coalition), its approximate solution is used for initial guesses to compute approximate solutions to all coalition structures with two coalitions.  Next, some Pareto optimal solutions are computed for all coalition structures, where any three-coalition coalition structure problem takes for initial guesses the solutions to the two-coalition problems for which it is a refinement. This is repeated until Pareto optimal solutions of the next finest coalition structure have been computed. If Stage 1 is performed, we refer to this as the \emph{bottom-up} approach, since we start at the lowest level of the coalition structure graph, and work upwards to the more refined levels. The rationale of approximating only some Pareto optimal solutions is to keep the computational cost limited, while providing a good initial guess to more efficiently find an accurate approximation of the full PF of the singleton coalition problem, which is the single most computationally complex part of the proposed algorithm (Stage 2). In Stage 3, the PS of the singleton coalition structure obtained in the previous stage is post-processed by checking, for all members of the PS of the grand coalitions, whether they are non-dominated also for the objective function corresponding to any given target coalition structure $CS$.

In the general case where~\eqref{eq:moo} has an objective function that cannot be expressed by Eq.~\eqref{eq:gen_obj_lin_agents} or is subject to constraints that  vary with the coalition structures (and not only vary with the actions of the coalitions), and when the full PF is of interest in decision-making, then numerical multi-objective optimization needs to be performed for all coalition structures. 
Within the loop over coalition structures, a MOO~\eqref{eq:moo}  needs to be numerically solved by repeatedly evaluating a physical solver describing the effects from agents operating on the subsurface resource of interest. This nesting property makes the problem potentially very computationally demanding, and it is crucial to both find an optimization method that does not require a very large number of model evaluations, and a physical model that is sufficiently simple without sacrificing too much accuracy. The algorithmic steps are summarized in Algorithm~\ref{alg:standard_MOO}, where we have indicated that prior information about the solution may be used to reduce the total computational cost. In contrast to the hierarchical approach in Algorithm~\ref{alg:hier_MOO}, obtaining such prior information (denoted $q^{CS}$ in Algorithm~\ref{alg:standard_MOO}) may be difficult in practice.

\begin{minipage}[t]{0.52\textwidth}
\begin{algorithm}[H]
\caption{\\ Hierarchical Multi-objective Optimization}
\label{alg:hier_MOO}
          \begin{algorithmic}[1]
          \phase{S1: Hierarchically compute some PS members.}
            \State Initialize design variables:\\ $q^{CS} \forall CS \in \Pi^{A}$
            \For{$i=1,\dots,|A|-1$}
                \For{$CS \in \Pi^A$ s.t. $|CS|=i$} 
                \State Partially solve:
                \State $q^{CS}~\leftarrow~\argmax_{q} (F_{C_1},\dots, F_{C_{|CS|}} | q^{CS} )$     
                \For{$CS' \leq CS$}
                		\State $q^{CS'} \leftarrow  q^{CS} \cup q^{CS'} $
                \EndFor 
                \EndFor 
                \EndFor
                \phase{S2: Compute PS of singleton CS.}
                \State $CS \leftarrow \{ \{a_1\},\dots, \{a_{|A|}\}\}$
                \State $q^{*} \leftarrow \argmax_{q} (F_{C_1},\dots, F_{C_{|CS|}} | q^{CS})$
                \State $\text{PS}_{CS} \leftarrow q^{*}$
                \phase{S3: Postprocessing to obtain all PS}
                \For{$CS \in \Pi^A$} 
                   \State $\text{PS}_{CS} \leftarrow \{q^{*} | q^{*} \text{ non-dominated in } CS \}$
                \EndFor
        \end{algorithmic}%
        
\end{algorithm}
\end{minipage}%
\hfill
\hspace{5pt}
\begin{minipage}[t]{0.48\textwidth}
\begin{algorithm}[H]
\caption{\\ Non-nested Multi-objective Optimization}
\label{alg:standard_MOO}
            \begin{algorithmic}[1]
                \For{$CS \in \Pi^A$} 
                \State Initialize design variables $q^{CS}$
                \State To desired accuracy, solve:
                \State $q^{*} \leftarrow \argmax_{q} (F_{C_1},\dots, F_{C_{|CS|}} | q^{CS})$
                \State $\text{PS}_{CS} \leftarrow q^{*}$
                \EndFor 
        \end{algorithmic}%
\end{algorithm}
\end{minipage}\\
\vspace{5pt}

Finally, we note that we deliberately have not touched upon the subject about how exactly the initial guesses or partial solutions should be used in numerical optimization. The proposed framework is not dependent on any particular type of optimization method, although it has been tacitly assumed that some kind of good initialization will indeed reduce the numerical cost.
The WSM  can be used to transform the MOO~\eqref{eq:moo} to a set of single-objective problems, but this strategy may become computationally very expensive~\cite{Pettersson_etal_24}. Alternatively, numerical methods to simultaneously approximate the full PF can be directly applied to~\eqref{eq:moo}. A very wide range of numerical MOO methods have been developed, and while a method should be selected based on the characteristics of the problem at hand for best performance, the exact choice of method is not the primary goal of this paper.

\section{Numerical results}
\label{sec:num_results}

We exemplify the MOO-based workflow for decision-making with cooperative games with two applications pertinent to subsurface systems: groundwater management and CO$_2$ storage. Specifically, we consider four test cases in total, each with three or four agents representing well operators, and between three and nine injection wells for fluid injection. In these examples, the coalitions seek to optimize their value, defined as the total amount of water or gas they inject. For simplicity, in the sense of limiting the variety of possible outcomes, the payoff of an agent is assumed to be equal to the amount that this agent injects in a given coalition structure. Note however that a different payoff by means of transferable utility (simply distributing the value differently from proportional to injection within a coalition) is possible to apply as a post-processing step of the existing numerical results.

As a first application, we consider an analytical superposition-of-wells model for groundwater flow where the Pareto optimal set is both convex and finite, and can efficiently be computed using the WSM (Eq.~\eqref{eq:WSM}) and linear programming~\cite{Marler_Arora_10}. Each agent is assumed to operate a single well, although the generalization to more than one well per agent does not introduce any conceptual changes to the model. 

In the second application, we present numerical results obtained with more efficient numerical methods compared to the same Bjarmeland formation test case shown in~\cite{Pettersson_etal_24} and highlight the variability in the solutions when using the selection criteria presented in Section~\ref{sec:pt_select}.
We also model the Bjarmeland formation with nine wells operated by three agents, different starting times for injection between wells, and longer injection times. 

\subsection{Groundwater resources}

We first rewrite the general multi-objective optimization problem~\eqref{eq:moo} in the special case of linearity, before applying it to a groundwater resource management problem. 
If the constrained multi-objective optimization problem~\eqref{eq:moo} is linear in the objective functions and constraints, the problem reduces to a multi-objective linear programming problem (MOLP), which provides interesting insights. To make this concrete and emphasize the relation between coalition structures and multi-objective optimization problems, let $M^{A} \in \mathbb{R}^{|A| \times N_{\text{dv}}}$ be the matrix such that the $i$th row of $M^{A}q$ denotes the objective function of agent $a_i\in A$, $i=1,\dots |A|$. Recall that $N_{\text{dv}}$ denotes the number of decision variables. For any coalition structure $CS \in \Pi^{A}$, let $M^{A2CS} \in \mathbb{R}^{|CS| \times |A|}$ be the matrix that maps from the agents $A$ to the current coalition structure $CS$, such that the $i$th row of $M^{A2CS} M^{A}q$ is the objective function of the $i$th coalition in $CS$. Typically, $[M^{A2CS}]_{ij} = 1$ if agent $a_{j}$ is in coalition $C_i \in CS$, and zero otherwise. For notational convenience, also set $M^{CS} \equiv M^{A2CS} M^{A}\in \mathbb{R}^{|CS| \times N_{\text{dv}}}$. With this notation, define the class of MOLPs for all $CS \in \Pi^A$ as
\begin{equation}
\label{eq:MOLP}
\max_{q\in\mathbb{R}^{N_{\text{dv}}}} M^{CS} q \quad
\text{s.t.}\quad B q \leq b \quad \text{and}\quad q_{\text{min}} \leq q \leq q_{\text{max}},
\end{equation}
where 
$B \in \mathbb{R}^{N_{\text{con}} \times N_{\text{dv}} }$, and $b \in \mathbb{R}^{N_{\text{con}}}$.
Here, and in the following, we assume that the optimization problems are well-posed~\cite{Dontchev_Zolezzi_06}, e.g., that the constraints admit the existence of a feasible solution. 
Applying the WSM~\eqref{eq:WSM} to the rows of $M^{CS}$ in~\eqref{eq:MOLP}, i.e., $c_{j}  \equiv \sum_{i=1}^{|CS|}\alpha_i M^{CS}_{i,j}$ for $j=1,2,\dots, N_{\text{dv}}$, results in the constrained single-objective linear programming (SOLP) problem,
\begin{equation}\label{eq:SOLP}
\max_{q\in\mathbb{R}^{N_{\text{dv}}}} c^{T}q \quad
\text{s.t.}\quad B q \leq b \quad\text{and}\quad q_{\text{min}} \leq q \leq q_{\text{max}},
\end{equation}
where $c = (c_1,\dots, c_{N_{\text{dv}}})^T$.

To apply this model to a groundwater resource management scenario, consider a homogeneous aquifer in 2D ($x,y$) as a groundwater flow model.  Water pumping or recharge can be performed independently by different agents, 
operating
$N_{\text{w}}$ water injection wells, located at $(x_k, y_k)$, $k=1,\dots, N_{\text{w}}$. 
Assume a water disposal scenario, where the injection rates are allowed to change only at the $N_{\text{t}}$ discrete times $\tilde{t} = (0, \Delta t, 2\Delta t,..., (N_{\text{t}}-1) \Delta t)$ for some uniform time window of length $\Delta t>0$. Let $q_{k,n}$ be the volumetric injection rate in well $k$ during the time interval $[(n-1)\Delta t, \ n\Delta t ]$, $n=1,\dots, N_{\text{t}}$. Furthermore, for any well $k$  we denote the injection rate change by
 \[
\Delta q_{k,n} \equiv \left\{
\begin{array}{ll}
q_{k, 1} & \mbox{if } n=1,\\
q_{k,n}-q_{k,n-1} & \mbox{if } n > 1
\end{array}
\right..
\]
Then we can express the change in hydraulic head using the Theis solution~\cite{Theis_35} as a superposition of $N_{\text{t}}$ independent injections (not necessarily positive) into the same well, each starting at a different time in $\tilde{t}$, and with injection ongoing at the time $t$ of evaluation, as follows:
\[
\Delta h(x,y,t) = 
\sum_{n=1}^{\ceil*{t/\Delta t}} \sum_{k=1}^{N_{\text{w}}} \frac{ \Delta q_{k,n}}{4 \pi T}W(\chi_{k,n}(x,y,t)),
\]
where $\ceil*{\cdot}$ denotes rounding to nearest larger integer, $T$ is transmissivity (area over time), $W(\chi)=\int_{\chi}^{\infty} \exp(-z)/z \textup{d}z$ is the well function for a confined aquifer, and the dimensionless group  $\chi_{k,n}$ is given by
\[
\chi_{k,n}(x,y,t) = \frac{S r_k^2(x,y)}{4 T (t-(n-1)\Delta t)},
\]
with dimensionless storage coefficient S, and radial squared distance $r_k^2=(x-x_k)^2+(y-y_k)^2$.
The injections are limited by maximum sustainable hydraulic head change (or pressure change; they are related via $\Delta p = g \rho_{\text{w}} \Delta h$), i.e., $\max_{x,y,t} \Delta h(x,y,t) \leq h_{\text{crit}}$ for some given $h_{\text{crit}}$. We note that the maximum hydraulic head change must occur at any of the discrete times $\tilde{t}$, and at a location on the radius of any of the injection wells. Choosing a single representative point on the radius of each well, the continuous hydraulic head constraint is reduced to $N_{\text{t}}N_{\text{w}}$ discrete constraints. In addition to uniform time discretization, we have assumed that the final time $t$ satisfies $t=N_{\text{t}} \Delta t$. The generalization when either of these assumptions does not hold is straightforward. 

The optimization problem can be written as a MOLP~\eqref{eq:MOLP} or SOLP~\eqref{eq:SOLP} where $M=\Delta t I_{N_{\text{w}}} \otimes (1,\dots, 1)$, the identity matrix of size $N_{\text{w}}$ is denoted by $I_{N_{\text{w}}}$, the symbol $\otimes$ denotes the Kronecker product, $c=\Delta t (1,1,\dots, 1)^T$, and $B=\tilde{B}(I_{N_{\text{w}}} \otimes D)$ with
\begin{multline}
\tilde{B}_{(i-1)N_{\text{w}} + j, (k-1)N_{\Delta t}+l} = \frac{Q_{\text{vol}}}{4\pi T}W_{ijkl}, \quad W_{ijkl} \equiv W(\chi_{k,l}(x_j, y_j,i\Delta t)),\\
\mbox{for } i=1,\dots, N_{\text{t}}, \quad j=1,\dots, N_{\text{w}}, \quad k=1,\dots, N_{\text{w}}, \quad l=1,\dots, i.
\end{multline}
Here, $D$ is the $(N_{\text{t}} \times N_{\text{t}})$-matrix with non-zeros only on the main and subdiagonals defined by $D_{i,j} = \delta_{i,j} - \delta_{i-1,j}$, $i,j=1,\dots, N_{\text{t}}$, where $\delta_{i,j}$ denotes the Kronecker delta. The factor $Q_{\text{vol}}$ is included to scale from injection rates in m$^3$/s to Mm$^3$/year. The problem can also be supplemented with constraints on all $q_{k,n}$, being bounded by the local constraints on injection and extraction  rates.

\subsubsection{Test case I: three agents}
For the first test case, we consider the superposition model with a spatial domain $[0, \ 10,000] \times [0, \ 10,000]$ [m$^2$] and set $S=1\times 10^{-5}$ [-], and $T=1\times 10^{-3}$ [m$^2$/s].  We assume three agents, each operating an injection well, located at positions $(2500, 2500)$, $(2500, 5000)$, and $(5000, 2500)$, respectively. Injection is performed simultaneously in all wells during a period of ten years, and each operator can change the injection rate annually, i.e., every $\Delta t =1$ year. The constraints on water disposal are $q_{\min}=40$ Mm$^3$/year and $q_{\max}=150$ Mm$^3$/year for all wells and all times, and $\Delta h \leq h_{\text{crit}} := 10,000$ m at all critical points in the domain. With this setup, the maximum constraint will never become active. Using the WSM~\eqref{eq:WSM},  results in a SOLP~\eqref{eq:SOLP} for each weight. Thanks to the convexity of the problem and modest computational cost, the full Pareto front is captured if a sufficiently fine discretization of the WSM weights is used. The linear programming problems are solved to a primal feasibility tolerance of $10^{-7}$ using the built-in primal-dual interior-point method in MATLAB R2021a~\cite{Mehrotra_92}, so the results presented are very accurate.

\begin{figure}[h]
    \centering  
{\includegraphics[width=1.00\textwidth]{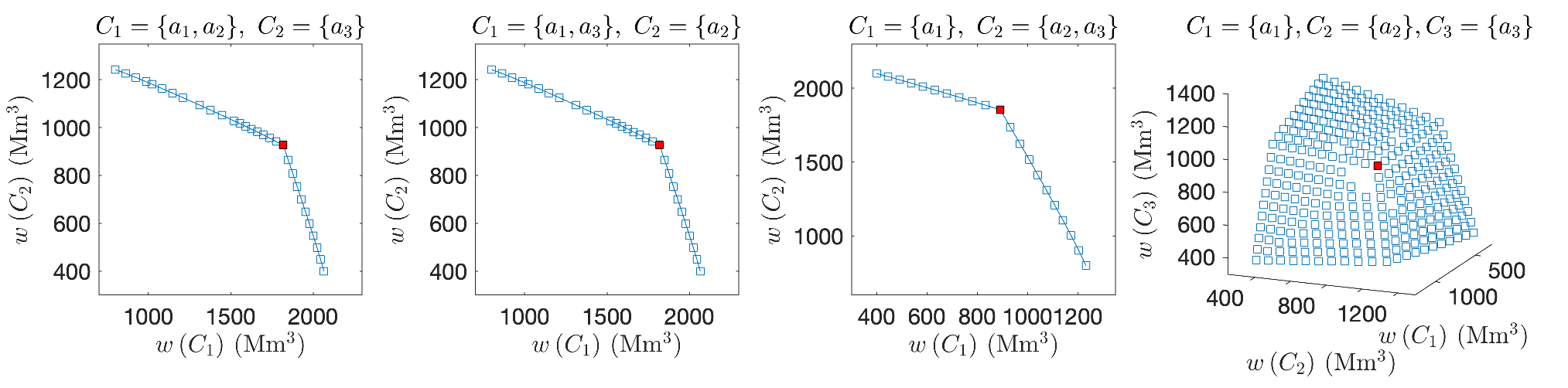}}
    \caption{Pareto fronts for the groundwater problem with three agents. The social welfare maximizing solution is indicated by a red marker. For ease of notation, we only indicate the the coalition $C$ of the value function, but the correct full notation is the embedded coalition $(C,CS)$.} %
    \label{fig:PF_groundwater}
\end{figure}

\begin{figure}[h]
    \centering  
{\includegraphics[width=1.00\textwidth]{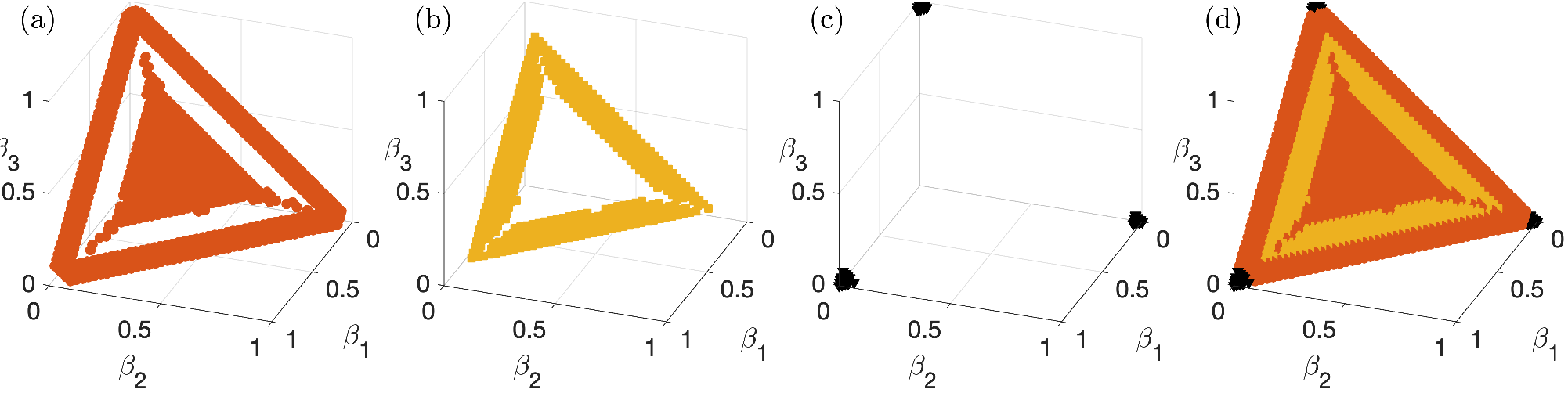}}
    \caption{Externalities for the groundwater problem with three agents. (a) Negative externalities; (b) mixed externalities; (c) zero externalities; (d) all externalities.} %
    \label{fig:extern_groundwater}
\end{figure}
Figure~\ref{fig:PF_groundwater} shows the Pareto fronts for the four coalition structures that have at least two coalitions. By virtue of Theorem~\ref{thm:hierarch_PS}, there exists a common point in all of the corresponding Pareto sets. This point on the Pareto fronts is social welfare maximizing and indicated by the red markers in Figure~\ref{fig:PF_groundwater}. It is also the unique solution for the grand coalition and the total injections are respectively 891, 927, and 927 Mm$^3$ for the three agents. The corresponding cooperative game has no externalities (i.e., it is a CFG) and all coalition structures are equally beneficial for all agents. The core consists of all coalition structures; there is no incentive to leave any coalition for any agent, since all coalition structures lead to identical injection schedules.

Next, we use the utopia point method described in Sect.~\ref{sec:pt_select} with $p=2$ and vary the weights $\tilde{\beta}_1,\tilde{\beta}_2,\tilde{\beta}_3$. Each set of weights corresponds to a single choice of points from all PFs, or, equivalently, a single game. Note that weights that are distinct but sufficiently close can result in the same game. We first characterize the games defined in this way by means of externalities, i.e., if and how the value of a coalition is affected by the structure of the remaining coalitions. For this test case, there are three externalities: 
\begin{equation}
\notag
\resizebox{1.0\linewidth}{!}{$
\begin{aligned}
\epsilon(\{ a_1\},\{ \{a_1 \}\{a_2,a_3\} \},\{ \{a_1 \}\{a_2\}\{a_3\} \}) &=
w(\{ a_1 \}, \{ \{ a_1 \}\{ a_2 , a_3 \} \}) -
w(\{ a_1 \}, \{ \{ a_1 \}\{ a_2 \}\{ a_3 \} \}) 
  \\
\epsilon(\{ a_2\},\{ \{a_2 \}\{a_1,a_3\} \},\{ \{a_1 \}\{a_2\}\{a_3\} \}) &= 
w(\{ a_2 \}, \{ \{ a_1, a_3 \}\{ a_2 \} \}) - w(\{ a_2 \}, \{ \{ a_1 \}\{ a_2 \}\{ a_3 \} \}) \\
\epsilon(\{ a_3\},\{ \{a_1, a_2 \}\{a_3\} \},\{ \{a_1 \}\{a_2\}\{a_3\} \}) &= 
w(\{ a_3 \}, \{ \{ a_1, a_2 \} \{ a_3 \} \}) 
- w(\{ a_3 \}, \{ \{ a_1 \}\{ a_2 \}\{ a_3 \} \})
\end{aligned}
$}
\end{equation}
Figure~\ref{fig:extern_groundwater} shows the externalities for different sets of positive weights with sums equal to 1 using Eq.~\eqref{eq:weighted_agents}.  Negative externalities mean that all three externalities are non-positive, and mixed externalities means that at least one externality is negative and at least one is positive. We observe that varying the weights yields a wide variety of situations, including negative, mixed, and zero externalities. Only positive externalities are not observed. Negative externalities are expected in games that model competition for finite resources, and it follows from Theorem~\ref{thm:hierarch_PS} that there exist solutions with no externalities. That we also observe mixed externalities, i.e., solutions where a coalition benefits from other coalitions collaborating to maximize their own value, may appear surprising and requires an explanation. The reason for this phenomenon is that while the utopia points are computed in the same way, they still vary between different coalition structures. Compared to when treated as a single coalition, when two subgroups of agents are treated as two distinct coalitions, their utopia points are more extreme since they each do not account for satisfaction of the goals of the other subcoalition. Hence, even if the merged coalition receives a larger weight than any of the individual coalitions, the penalty in the selection criterion~\eqref{eq:utopia_norm} can be relatively smaller due to the smaller distance to the joint utopia point. The result is that a third external coalition can receive a higher value of its partition function when the two sub-coalitions in question are merged into a single coalition compared to when they are competing disjoint coalitions. As an example, the case with agent weights $(\tilde{\beta}_{a_1}, \tilde{\beta}_{a_2}, \tilde{\beta}_{a_3})=(0.45, 0.45, 0.1)$ shown in Table~\ref{table:pfg_3agents} displays the following mixed externalities:
\begin{equation}
\notag
\begin{aligned}
\epsilon(\{ a_1\},\{ \{a_1 \}\{a_2,a_3\} \},\{ \{a_1 \}\{a_2\}\{a_3\} \}) &= 891 - 1025 < 0, \\
\epsilon(\{ a_2\},\{ \{a_2 \}\{a_1,a_3\} \},\{ \{a_1 \}\{a_2\}\{a_3\} \}) &= 927 - 1040 < 0,\\
\epsilon(\{ a_3\},\{ \{a_1, a_2 \}\{a_3\} \},\{ \{a_1 \}\{a_2\}\{a_3\} \}) &= 449-400 > 0.
\end{aligned}
\end{equation}

\begin{table}[h]
\centering
{\renewcommand{\arraystretch}{0.5}
\begin{tabular}{lrrr}
\multicolumn{1}{c}{} & \multicolumn{3}{c}
{$\tilde{\beta}$} \\
\multicolumn{1}{c}{} & $[1/3, 1/3, 1/3]$ & $[0.45,  0.45, 0.1]$ & $[0.98,\ 0.01, \ 0.01]$ \\
\multicolumn{1}{c}{} & (neg. ext.) & (mix. ext.) &  (no ext.) \\
\midrule
$w(\{ 1, 2, 3 \}, \{\{1, 2, 3\} \})$ & 2745 
$\left\{ \begin{tabular}{c@{}} 891 \\
  927 \\
  927 
  \end{tabular}
  \right.$
  & 2745 
$\left\{ \begin{tabular}{c@{}}
891\\
927\\
927
    \end{tabular}
  \right.$ &
  2745 
  $\left\{ \begin{tabular}{c@{}} 891 \\
  927 \\
  927 
  \end{tabular}
  \right.$
  \\
\midrule
$w(\{ 1, 2\}, \{\{1,2\},\{3\} \})$ &  1849 
$\left\{
\begin{tabular}{c@{}}
908\\
941
  \end{tabular}
  \right.$
& 2045 
$\left\{
\begin{tabular}{c@{}}
  1014\\
  1031
  \end{tabular}
  \right.$
& 2065 
$\left\{
\begin{tabular}{c@{}}
 1025 \\
1040
  \end{tabular}
  \right.$
\\
$w(\{ 3\}, \{\{1, 2\},\{ 3\} \})$ & 865 & 449  & 400  \\
\midrule
$w(\{ 1, 3\}, \{\{1, 3\},\{2\} \})$ & 1849 
$\left\{
\begin{tabular}{c@{}}
  908 \\
  941
  \end{tabular}
  \right.$
& 1818  
$\left\{
\begin{tabular}{c@{}}
891\\
927
  \end{tabular}
  \right.$
& 2065  
$\left\{
\begin{tabular}{c@{}}
   1025\\
   1040
\end{tabular}
  \right.$
\\
$w(\{ 2\}, \{\{1, 3\},\{2\} \})$ & 865 & 927 & 400 \\
\midrule
$w(\{ 1\}, \{\{1\},\{2 ,3\} \})$ &  833 & 891 & 1232\\
$w(\{ 2, 3\}, \{\{1\},\{2, 3\} \})$ &  1884 
 $\left\{
\begin{tabular}{c@{}} 
  942\\
  942
   \end{tabular}
  \right.$ 
& 1854  
  $\left\{
\begin{tabular}{c@{}}  
  927\\
  927
   \end{tabular}
  \right.$   
& 800  
$\left\{
\begin{tabular}{c@{}} 
 400\\
 400 
\end{tabular}
  \right.$ 
\\
\midrule
$w(\{ 1\}, \{\{1\},\{2\},\{3\} \})$ &  891 &  1025 & 1232  \\
$w(\{ 2\}, \{\{1\},\{2\},\{3\} \})$ & 927 &  1040 & 400 \\
$w(\{ 3\}, \{\{1\},\{2\},\{3\} \})$ & 927 &  400 & 400 \\
\midrule
Welfare-maximizing CS   & \multicolumn{1}{c}{ 
\begin{tabular}{@{}c@{}} $\{\{1, 2, 3\} \}$ \\ 
$\{\{1\},\{2\},\{3\} \}$
\end{tabular}
}
&
\multicolumn{1}{c}{
\begin{tabular}{c@{}} $\{\{1, 2, 3\} \}$ \\ $\{\{1, 3\},\{2\} \}$\\
$\{\{1\},\{2 ,3\} \}$
\end{tabular}
}
&
\multicolumn{1}{c}{
$\{\{1, 2, 3\} \}$
}
\\
\bottomrule
\end{tabular}
\caption{Payoffs and coalition values (total water disposal in Mm$^3$) of selected three-agents games with different weights assigned to the agents, and coalition values $w((C,CS))=F_C^{\beta}$ from Eq.~\eqref{eq:utopia_norm}.}
\label{table:pfg_3agents}
}
\end{table}
Table~\ref{table:pfg_3agents} displays the coalition values for all coalition structures for three representative games  with different externalities, and weights assigned to the agents as in Figure~\ref{fig:extern_groundwater}. Next, we briefly analyze the outcomes of the three games. For the first game with negative externalities, the grand coalition is not stable since any two agents could decide to leave the coalition and create their own coalition with improved payoff for both, at the expense of the remaining agent. No agent benefits from leaving a two-agent coalition. Thus, the core, consisting of all stable coalition structures, is the set of all coalition structures that contain a two-agent coalition. Interestingly, the social welfare of these coalitions is smaller than for the two remaining unstable coalition structures. This is in contrast to the situation for CFGs, where a necessary condition for a coalition structure to belong to the core is that it is social welfare maximizing~\cite[Prop.~2.2.1]{Chalkiadakis_etal_22}.

For the game with mixed externalities in Table~\ref{table:pfg_3agents}, the only stable coalition structure is the one with only singleton coalitions. The small weight on agent $a_3$ leads to the minimum guaranteed injection of 400 Mm$^3$. A smaller lower injection constraint would probably lead to lower injected volumes for this agent. From a social welfare point of view, this coalition structure only results in about 90 \% of the value that could be attained in an ideal scenario.

\subsubsection{Test case II: four agents}
\begin{figure}[h]
    \centering  
{\includegraphics[width=1.00\textwidth]{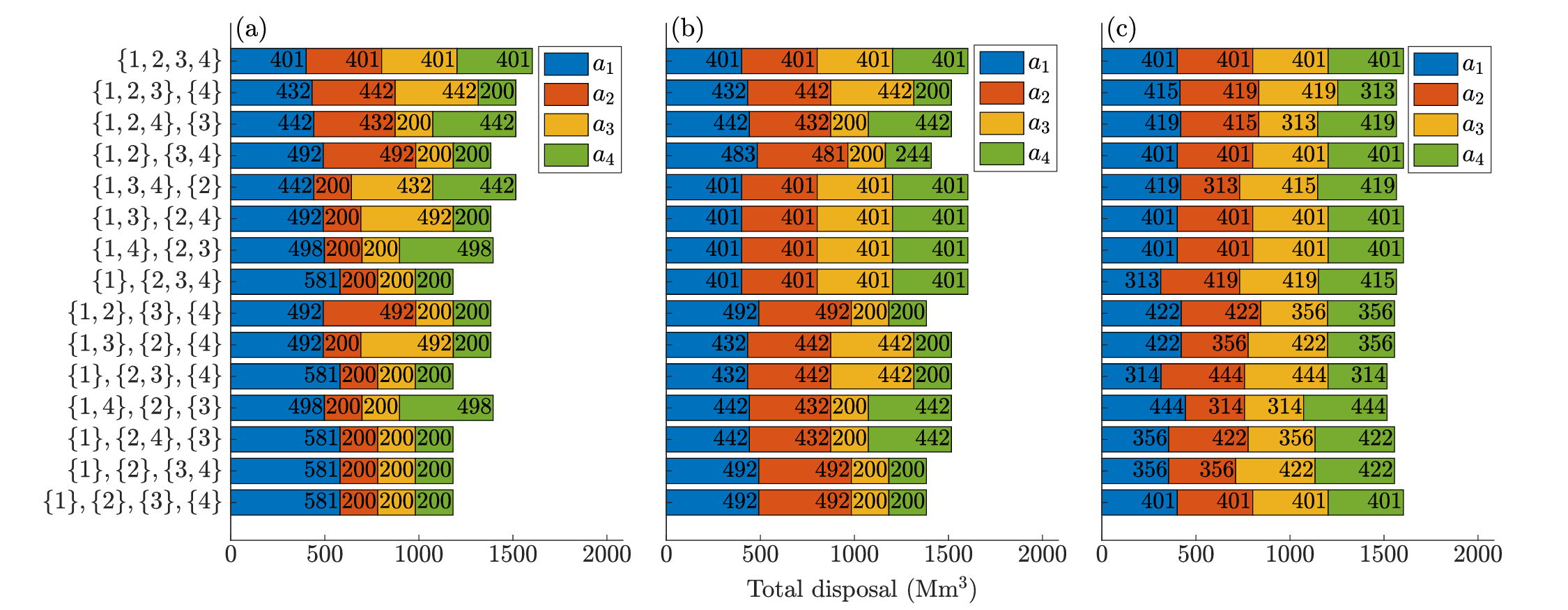}}

    \caption{Total injection  volumes for the weighted utopia point selection criterion, Eq.~\eqref{eq:utopia_norm}, with coalition weights~\eqref{eq:weighted_agents} given by (a) $\tilde{\beta} = [0.97, 0.01, 0.01, 0.01]$, (b) $\tilde{\beta} =[0.45, 0.45, 0.05, 0.05]$, and (c) $\tilde{\beta} = [0.25, 0.25, 0.25, 0.25]$.} %
    \label{fig:bar_Theis_4wells}
\end{figure}

Finally, we consider the groundwater problem with four wells with locations $(2500, 2500)$, $(2500, 5000)$, $(5000, 2500)$, and $(5000, 5000)$, each operated by an independent agent. Injection is performed during a period of 5 years, and the injection rates can change every year. Otherwise, the setup is similar to the first test case. The PFs are similar to those of the first test case, and not included here. To visualize representative games for the 15 coalition structures, we present results for selected agent weights using the Pareto front selection criterion~\eqref{eq:utopia_norm} in Figure~\ref{fig:bar_Theis_4wells}. This test case is completely symmetric in the four wells. In Figure~\ref{fig:bar_Theis_4wells}~(a), the coalition that contains $a_1$ is always favored due to the weight $\tilde{\beta}_{a_1}=0.97$.
Putting the same set of weights to different agents gives identical results after re-labeling of the coalitions due to the symmetry of the setup. Hence, the presented results are representative for assigning most weight to a single agent (a), assigning strong weight to two agents (b), and assigning equal weight to all agents (c). 

Due to the symmetry of the well locations, the (welfare maximizing) grand coalition results in equal payoffs for all agents of about 400 Mm$^3$ per agent. When a single agent is favored and receives a payoff of 581 Mm$^3$, all other agents achieve the minimum of 200 Mm$^3$, which is identical to the minimum constraint ($q_{\min} = 40$ Mm$^3$/year during 5 years). This suggests that by allowing a lower constraint, the outcome would be more extreme in the sense of allowing higher payoff of the favored coalition at the expense of the others.

Concluding the groundwater examples, we observe that, despite the physical problem being linear in both the objective function and the constraints,
we get a surprisingly complex class of games, displaying different characteristics in terms of externalities. In addition to these qualitative differences, which impacts the choice of coalition structure generation algorithm for many-agent problems, we also observe significant quantitative variability in the payoffs, depending on the game as a result of Pareto front selection criterion. 

\subsection{Large-scale CO$_2$ storage in Barents Sea}
Next, we present numerical results for a large-scale prospective CO$_2$ storage site in Barents Sea, the Bjarmeland formation, where injection is limited by pressure buildup.
It is sufficient to monitor pressure during the injection phase only since the maximum pressure in the physical model will be attained during this phase. While the long-term migration of the CO$_2$ plume is important, it has no impact on the objective functions in this work and will therefore not be further investigated. A fully implicit vertically integrated black-oil type model with two phases (CO$_2$ and brine) is used to compute the migration of CO$_2$ in the reservoir~\cite{Nilsen_etal_16}. Appropriate simplifications are employed, including coarse grids and uniform rock properties, as described in~\cite{Allen_etal_17}. For more details of the physical considerations and numerical setup, we refer to~\cite{Allen_etal_17} and the open-source test cases in MRST, c.f.\ \url{https://www.sintef.no/projectweb/mrst/modules/co2lab/}, that are used for all numerical experiments.

\subsubsection{Test case III: Bjarmeland formation with three wells}
\label{sec:Bjarmelans_3wells}
\begin{figure}[h]
    \centering  
{\includegraphics[width=0.48\textwidth]{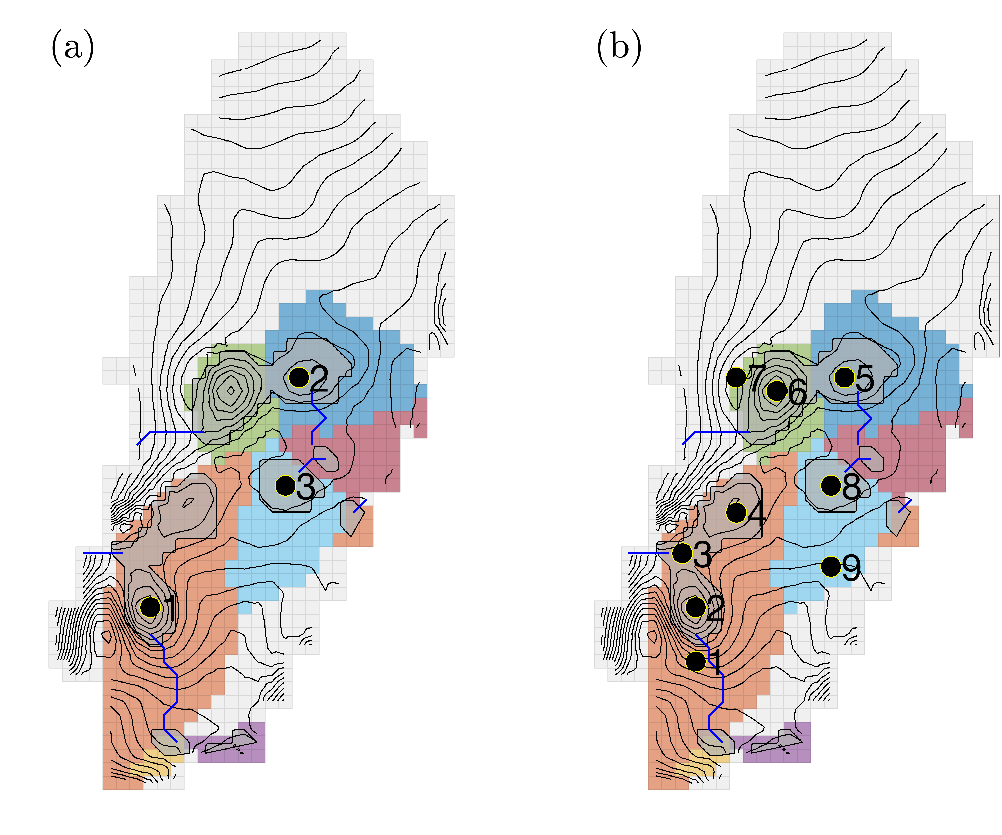}}
    \caption{Well locations for the Bjarmeland problem, (a) with three agents and three wells, and (b) three agents and nine wells. The colored areas are the catchment regions.} %
    \label{fig:Bjarmeland_well_pos}
\end{figure}

Three injection wells, each operated by an independent agent, are located at the peaks of some of the largest structural traps of the Bjarmeland formation, as shown in Figure~\ref{fig:Bjarmeland_well_pos} (a).  The well coordinates are identical to three of the four wells investigated in~\cite{Allen_etal_17}, where the omitted well only allowed limited injection. For every coalition structure, we seek the injection rates that maximize the total amount of CO$_2$ for every coalition. For simplicity, we assume for this first Bjarmeland test case that all wells start injection simultaneously, inject for 15 years,  and that they can change their injection rates at the same predefined times every three years. With five injection intervals per well, there is a total of 15 decision variables (one per well and injection interval) for each optimization problem.  
 All wells are subject to the same constraints imposed due to supply of CO$_2$ and minimum injection rates in the current numerical setup, but there is no restriction on varying the constraints between the wells in the proposed framework. We assume  minimum and maximum injection rates of respectively 0.24 and 7 Mton/year, where the former is supposed to represent a constraint based on economic feasibility and the latter a supply constraint.  The maximum rates are intentionally chosen to be less restrictive than the physical pressure buildup constraint, to be defined next. The problem would otherwise simplify so that the physical model becomes redundant, in the sense that it has no effect on the numerical results. 
The nonlinear part of the constraint function $g$ represents an upper limit on the reservoir pressure $p_{\text{res}}$ set to 90~\% of the overburden pressure $p_{\text{ob}}$, evaluated pointwise in space:
 \[
 g_i = \frac{p_{ \text{res}}(\mathbf{x}_{i}) }{p_{\text{ob}}(\mathbf{x}_{i}) } - 0.9, \quad \mbox{for } i \mbox{ so that } \mathbf{x}_{i} \in \mathbf{X}_{\text{grid}},
 \]
 where $\mathbf{X}_{\text{grid}}$ denotes the discrete spatial grid of the numerical reservoir model. This constraint was previously also imposed in optimization of the Bjarmeland formation in~\cite{Allen_etal_17, Pettersson_etal_24}.

 The MOO problem~\eqref{eq:moo} for each of the coalition structures is solved with the nondominated sorting algorithm NSGA-II, that uses a crowding distance measure for maintaining population diversity, and elitism by keeping track of identified non-dominated solutions for improved performance. This is a robust and popular algorithm that has been widely used since its introduction in~\cite{Deb_etal_02}. 
NSGA-II displays good diversity properties, in the sense that the members of the PF display relatively uniform spacing. Finding the extrema (end points) of the PF with NSGA-II may however require a very large number of function evaluations. As a remedy, we first find the extrema by solving a single-objective optimization problem for each coalition, as follows~\cite{Rostamian_etal_25}. By Lemma~\ref{lemma:WSM_PF}, we can find a single point on the desired PF by the WSM~\eqref{eq:WSM} and using any positive weights. To approximately determine the extreme points, we set the WSM weights $\alpha$ to a value close to 1 for the target coalition, and some small number $\epsilon$ for all other coalitions. In all numerical experiments, we use $\epsilon=0.001$. Putting essentially all the weight on the target coalitions implies that we should obtain the Pareto solution that is the most beneficial for the target coalition. Non-zero weight on remaining coalitions ensures that the solutions are indeed on the PF so that the non-target coalition do not just get assigned any feasible value that could be improved upon without detriment to the target coalition. The single-objective optimization problem that results from an application of the WSM is then solved using competitive swarm optimization (CSO), a variant of particle swarm methods with enhanced population diversity where the superior half of the population  is transferred directly to the next generation, and the inferior half is updated based on the superior individuals~\cite{Cheng_Jin_14}. 
 
 To reduce the number of model evaluations until numerical convergence, the initial population consists of Latin Hypercube samples over the entire decision variable space, complemented by a population member given by the minimum values of all decision variables. The latter ensures that at least one member of the initial population is a feasible solution (i.e., that does not violate any constraint). For both the constrained MOO and the single-objective optimization, we employ NSGA-II and CSO implementations from the MATLAB Platform for Evolutionary Multi-objective Optimization (PlatEMO)~\cite{Tian_etal_17}.

The MOOs corresponding to the Bjarmeland problem setup satisfies the assumptions of Theorem~\ref{thm:hierarch_PS}, so the Pareto sets satisfy the hierarchical subset property that should be utilized for reduced numerical cost in the approximation of the solutions. To investigate the relative performance of the hierarchical MOO methods proposed in Section~\ref{sec:algorithms}, we will compare results for the non-nested MOO Algorithm~\ref{alg:standard_MOO} with independent computation of the PFs, to results where we use, respectively, the bottom-up and top-down versions of Algorithm~\ref{alg:hier_MOO}. The same combination of NSGA-II and CSO described above, is used for each individual MOO, whether the hierarchical properties of the PSs have been used or not.
 
\begin{figure}[h]
    \centering  
{\includegraphics[width=1.00\textwidth]%
{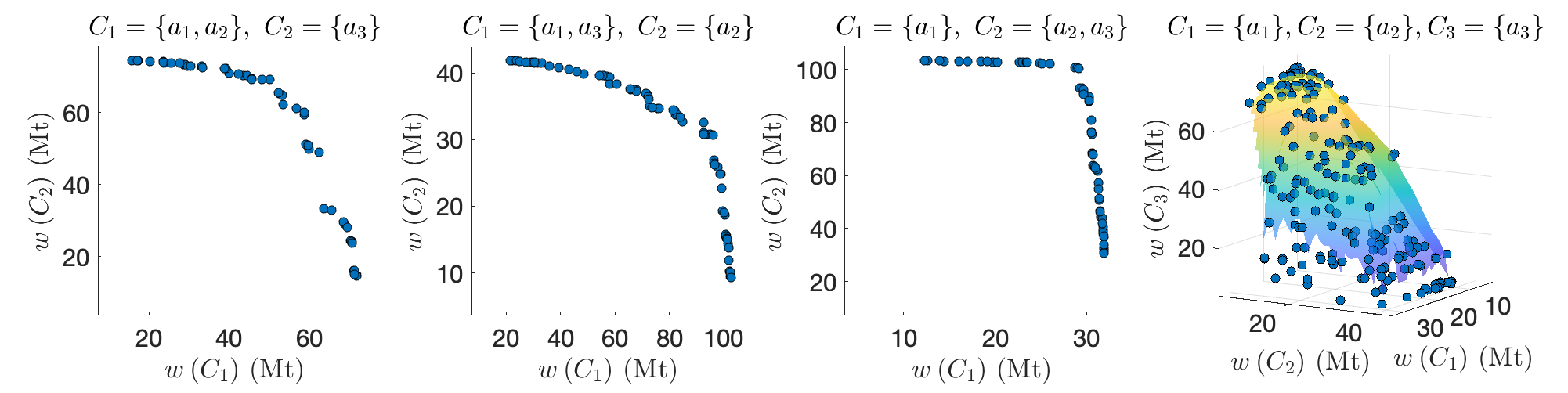}}
    \caption{Pareto fronts for the Bjarmeland problem with three agents and three wells. Independent simulations, a total of 111,000 model evaluations.} %
    \label{fig:PF_Bjarmeland_3wells}
\end{figure}

The PFs computed independently (Algorithm~\ref{alg:standard_MOO}) for this test case are shown in Figure~\ref{fig:PF_Bjarmeland_3wells}.

 \begin{figure}[h]
    \centering  
{\includegraphics[width=1.00\textwidth]
{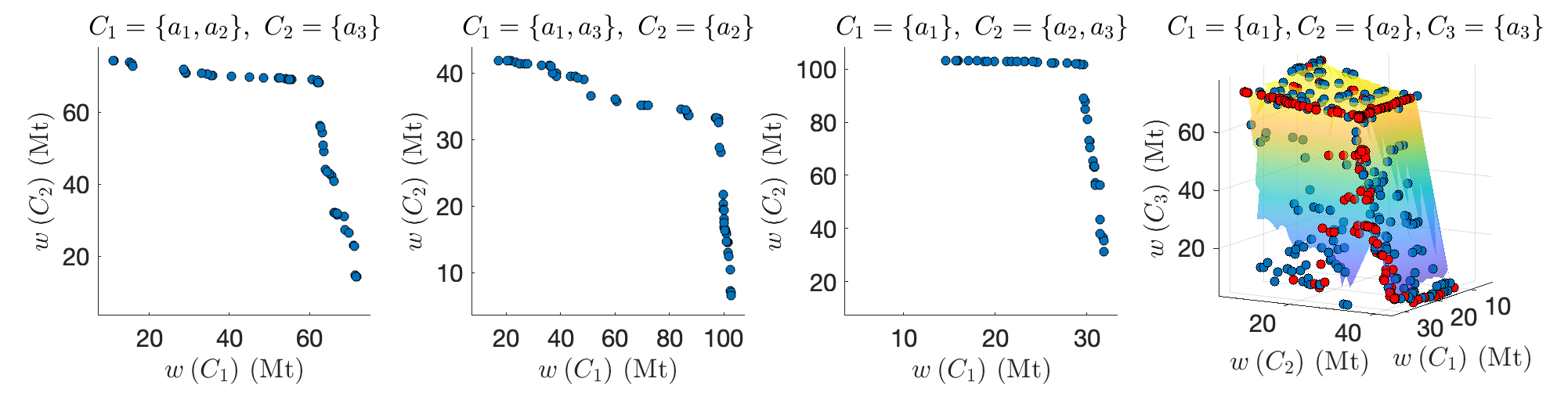}}
    \caption{Pareto fronts for the Bjarmeland problem with three agents and three wells. Hierarchical simulations, bottom-up approach, a total of 91,000 model evaluations.} %
    \label{fig:PF_Bjarmeland_3wells_bottom_up}
\end{figure}
\begin{figure}[h]
    \centering  
{\includegraphics[width=1.00\textwidth]
{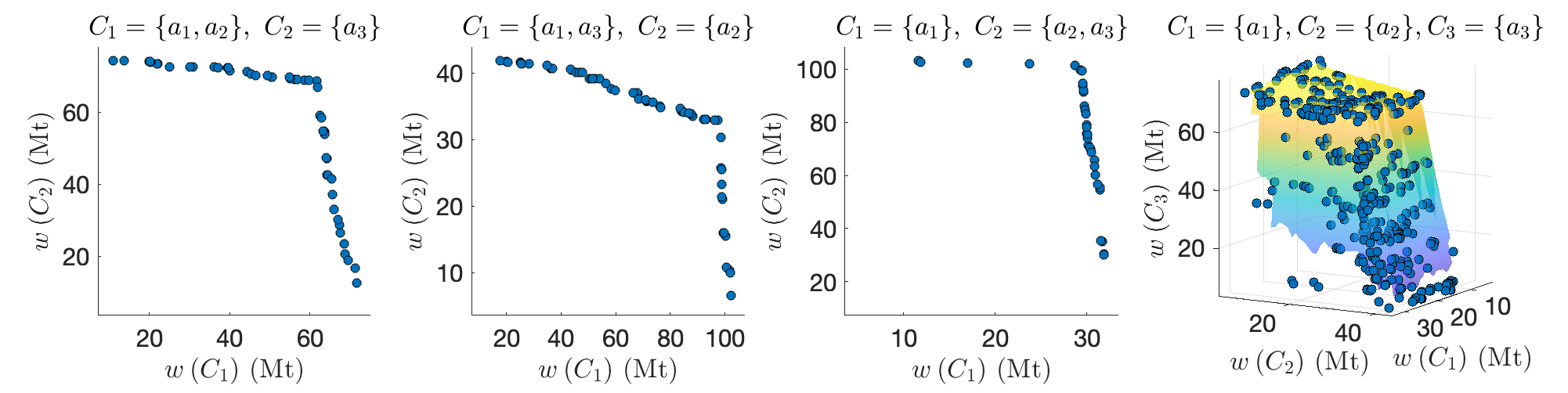}}
    \caption{Pareto fronts for the Bjarmeland problem with three agents and three wells. Hierarchical simulations, top-down approach, a total of 60,000 model evaluations.} %
    \label{fig:PF_Bjarmeland_3wells_top_down}
\end{figure}

A surface has been fitted to the rightmost PF for improved visibility of its shape. It highlights that some regions are sparsely populated. Whether these regions are due to an insufficient number of computed solutions (i.e., too small population size in the evolutionary algorithm), or actual reflections of the true PF, is not immediately clear.

 Next, we use the bottom-up approach described in Section~\ref{sec:algorithms} with Stage 1 included in Algorithm~\ref{alg:hier_MOO}. We start with the PFs of lowest dimension, corresponding to the coarsest coalition structure, and use the numerical results in the computation of higher-dimensional PFs corresponding to more refined coalition structures. By Theorem~\ref{thm:hierarch_PS}, the single-objective optimization problem for the grand coalition yields an approximation of the welfare-maximizing solution for all PFs. This solution, and the relevant end-point solutions computed using CSO, are added to the initial sets of the NSGA-II to hierarchically compute the PFs, finishing with the 3D PF approximation of the singleton coalition structure. The PFs are shown in Figure~\ref{fig:PF_Bjarmeland_3wells_bottom_up}. Blue markers indicate that the solutions have been computed for the current coalition structure, and red markers (rightmost PF) that they have been computed for a coarser coalition structure.

In the application of the top-down approach described in Section~\ref{sec:algorithms} and Algorithm~\ref{alg:hier_MOO}, we directly solve the MOO corresponding to the singleton coalition structure. The initial decision variable set is complemented with the same maximized social-welfare and end-point approximations used in the bottom-up approach. All lower-dimensional PF approximations are performed by post-processing of the singleton coalition PF, without any additional numerical simulations. The results are shown in Figure~\ref{fig:PF_Bjarmeland_3wells_top_down}.

Assessment of the performance of the methods employed to generate the PF approximations requires both consideration of the quality of the solution, and the computational time. The latter is to a very good approximation linear in the number of model evaluations, and hence straightforward to assess. The quality of the PF is more challenging to estimate, as there is no unique measure of quality that captures all aspects of desired PF approximation properties~\cite{Zitzler_etal_08}. Furthermore, evolutionary methods employed in this work, are stochastic and hence call for statistical PF metrics. One relatively simple and widely used metric for PF approximation is the hypervolume indictor, estimating the measure of the objective space enclosed by the generated PF and a user-defined reference point~\cite{Zitzler_Thiele_98}. In this work, we have compared the hypervolume indicators of the three sets of four PF approximations shown in Figures~\ref{fig:PF_Bjarmeland_3wells}--\ref{fig:PF_Bjarmeland_3wells_top_down}. According to this metric, the 2D PFs are all similar for all three methods, or slightly better with the two hierarchical methods. The hypervolumes corresponding to the 3D PFs are 10 \% larger for both hierarchical methods, compared to the non-nested method. The two hierarchical approaches are similar in terms of hypervolume indicator metrics. Since the PF approximation using the top-down approach was generated with the smallest number of function evaluations, we conclude that this is the most efficient method for the Bjarmeland problem. Hence, we will henceforth use the corresponding PFs when investigating the cooperative games.

\begin{figure}[h]
    \centering  
{\includegraphics[width=0.96\textwidth]{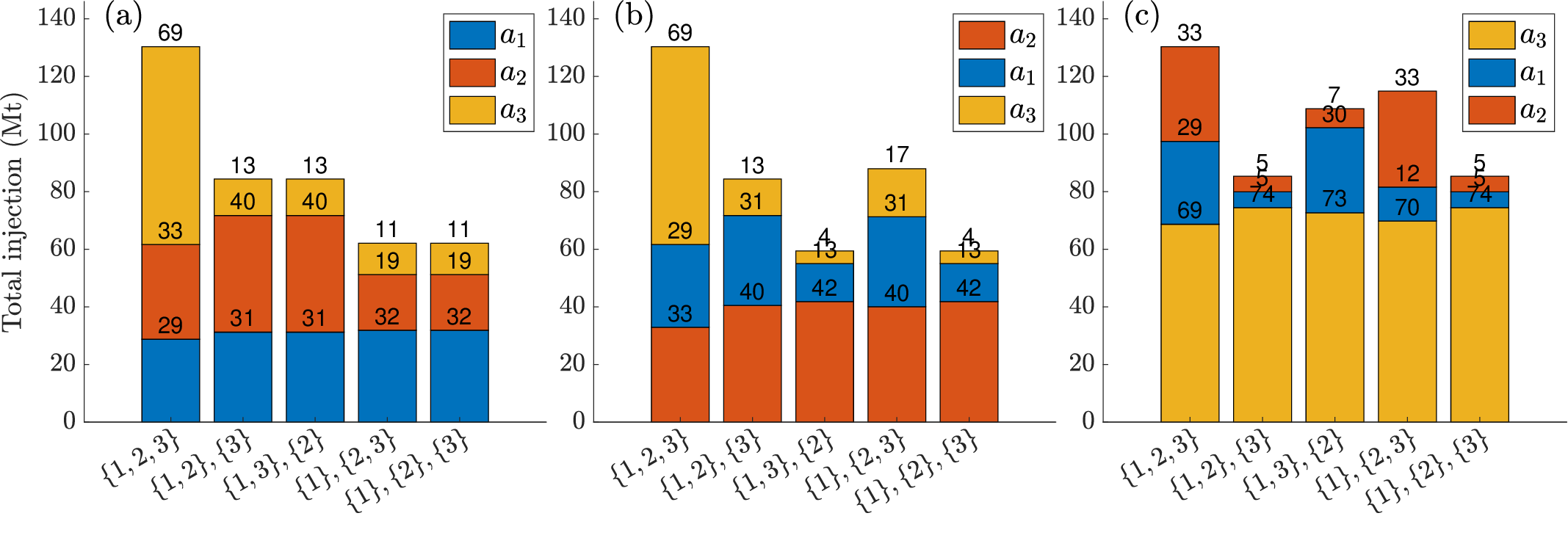}}
    \caption{Total injection  when consistently choosing the Pareto optimal solution that maximizes the value of (a) agent $a_1$, (b) agent $a_2$,  and (c) agent $a_3$.} %
    \label{fig:Bj_3w_rates_agent_ops}
\end{figure}

\begin{figure}[h]
    \centering  
{\includegraphics[width=1.00\textwidth]{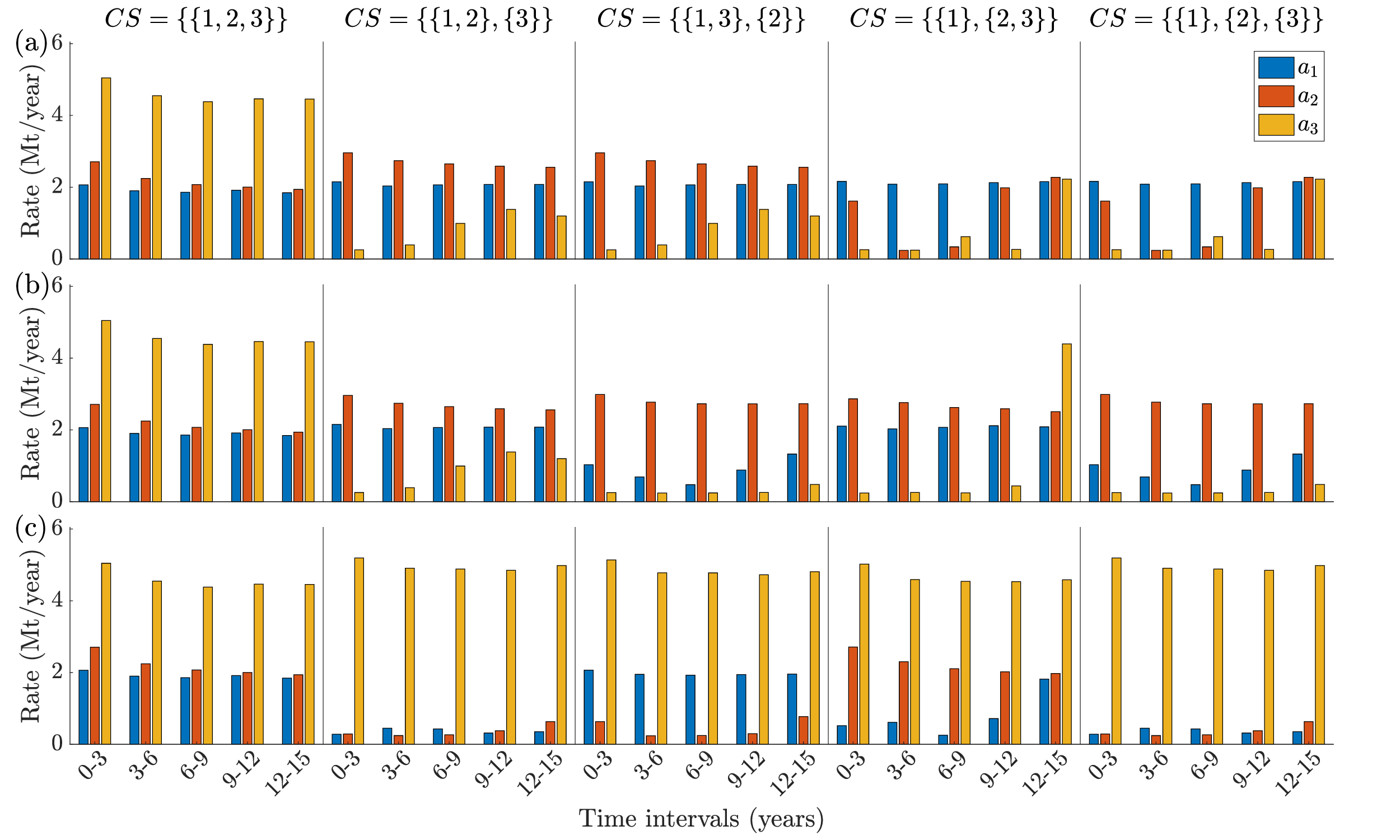}}
    \caption{Injection rates when consistently choosing the Pareto optimal solution that maximizes the value of an individual agent. (a) Maximizing for $a_1$, (b) Maximizing for $a_2$, (c) Maximizing for $a_3$.} %
    \label{fig:Bj3w_annual_rates}
\end{figure}

We assemble the cooperative games by applying the PF selection criterion~\eqref{eq:favor_ind_agent} to the PFs in Figure~\ref{fig:PF_Bjarmeland_3wells_top_down} to systematically favor agent $a_1$, $a_2$, and $a_3$, respectively. The payoffs, assumed equal to the individual total injections, are shown in 
Figure~\ref{fig:Bj_3w_rates_agent_ops}. By construction, the grand coalition with its single optimum is unaffected by the selection criterion and provides the total maximum attainable. The favored agent always benefits from not joining any coalition, but the gain is very small compared to the overall loss of total injection volumes. Less than 50 \% of the maximum possible amount of CO$_2$ is injected when either $a_1$ or $a_2$ is favored, and the outcome is that this agent chooses not to collaborate. That number is 65 \% when $a_3$ is favored, which is also far from maximal utilization of storage resources. In all PF candidate solutions, only the physical pressure constraints are active. This suggests that there would be no gain for any given agent if another agent had a smaller minimum injection rate. Finally, Figure~\ref{fig:Bj3w_annual_rates} shows the annual injection rates for the optimized solutions. They are overall relatively even over time, but there are some exceptions, in particular for agent $a_3$ where dramatic changes are observed. Some regularization of the rates could be investigated as it may be the case that more even rates would not significantly deteriorate the overall performance.

\subsubsection{Test case IV: Bjarmeland formation with nine wells}

We now consider the Bjarmeland formation, with nine wells operated by three agents and located as shown in Figure~\ref{fig:Bjarmeland_well_pos} (b). We consider the schedule, where wells 2, 3, and 5-9 start injection at the same time and inject for 40 years. Wells 1 and 4 start injecting five years after the others, and inject for a total of 35 years. The injection rates can change every 5 years. Agent 1 operates wells 1-4, agent 2 operates wells 5-7, and agent 3 wells 8-9. Each MOO has 70 decision variables (one per well and injection well) and is clearly more complex than the previous test case. The constraints on pressure and minimum/maximum annual injection rates remain the same.
Injected CO$_2$ that migrates outside the catchment regions (colored regions in Figure~\ref{fig:Bjarmeland_well_pos}) will eventually leak out~\cite{Allen_etal_17}. Nevertheless, in this work, we focus on the pressure constraints only and do not track the migration of the CO$_2$ plume. 

\begin{figure}[h]
    \centering  
{\includegraphics[width=1.00\textwidth]{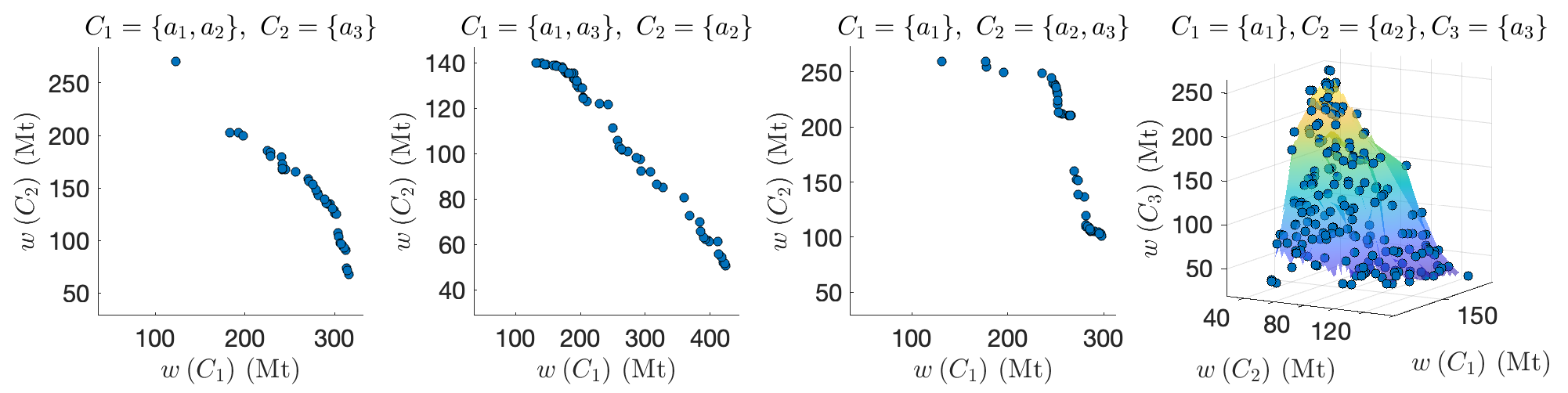}}
    \caption{Pareto fronts for all coalition structures with at least two coalitions for the Bjarmeland test case with 9 wells. Independent simulations, a total of 100,000 model evaluations. The axes intercepts coincide with the minimum injection constraints.} %
    \label{fig:Pareto_BF_Bj9}
\end{figure}

\begin{figure}[h]
    \centering  
{\includegraphics[width=1.00\textwidth]{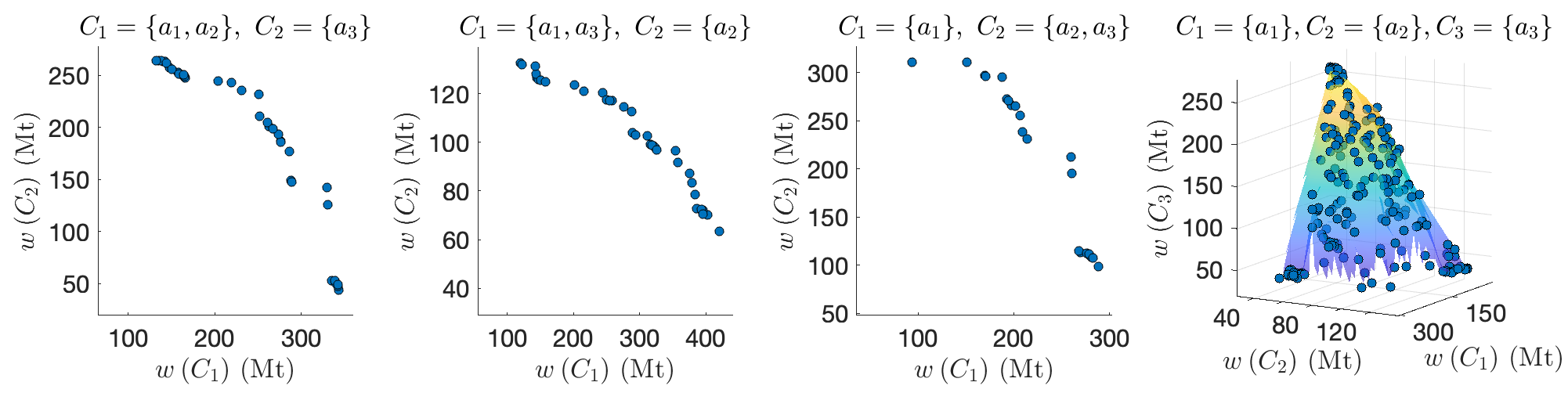}}
    \caption{PFs for all coalition structures with at least two coalitions for the Bjarmeland test case with 9 wells. Hierarchical simulations, top-down approach with a total of 70,000 model evaluations. The axes intercepts coincide with the minimum injection constraints.} %
    \label{fig:Pareto_top_down_Bj9}
\end{figure}

The PFs of the four coalition structures with at least two coalitions computed independently are shown in Figure~\ref{fig:Pareto_BF_Bj9}. The same PFs with the hierarchical top-down approach are displayed in Figure~\ref{fig:Pareto_top_down_Bj9}.
For both the non-nested and the hierarchical top-down methods, 5000 physical model evaluations with CSO are used to find the extreme points of the PFs (population size 50, and 100 iterations). Then, 15,000 (population size 200, 75 iterations) model evaluations are employed to find the PF with NSGA-II, for each coalition structure with the non-nested approach, resulting in a total of 100,000 model evaluations. For the top-down approach, 20,000 (population size 200, 100 iterations) model evaluations are used for the singleton coalition structure only, in total 70,000 model evaluations.
Comparing the quality of the two sets of PFs is not straightforward, as described previously in Section~\ref{sec:Bjarmelans_3wells}, although the hierarchical approach seems to capture larger parts of the PFs. This is confirmed by the hypervolume indicator, according to which the 2D fronts are at least as good for the top-down approach as for the non-nested approach, and the 3D PF approximation encloses a hypervolume that is about 50 \% larger for the top-down approch. Different quality metrics, and different design choices, e.g., population size, or using an altogether different MOO method, may yield a more accurate solution.

\begin{figure}[h]
    \centering  
{\includegraphics[width=1.00\textwidth]{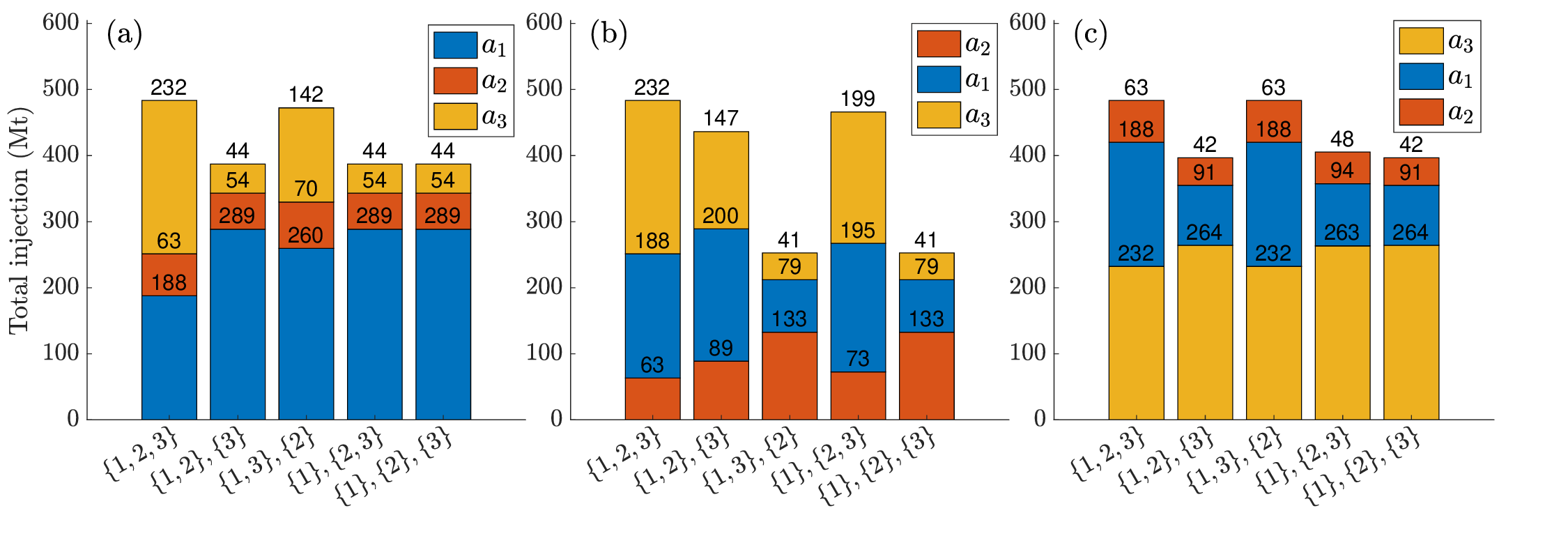}}

    \caption{Total injection  when consistently choosing a Pareto optimal solution that maximizes the value of (a) agent $a_1$, (b) agent $a_2$, and (c) agent $a_3$. 
    } %
    \label{fig:tot_inj_bars_9wells}
\end{figure}

Selected games from the PFs in Figure~\ref{fig:Pareto_top_down_Bj9} are shown in Figure~\ref{fig:tot_inj_bars_9wells}. Compared to the previous test case with three wells, there are now nine wells, so a direct comparison is not meaningful. Depending on which agent is favored, that agent's payoff is maximized at the expense of the overall degree of utilization of the subsurface resource as a whole, defined as the maximum total possible amount injected and realized by the grand coalition. In particular, $a_2$ can make a significant gain but at the loss of total subsurface storage efficiency, which is then not much more than 50 \% of the storage potential.

\section{Discussion}
The combination of cooperative game models with value functions determined by the solutions to multi-objective optimization problems subject to complex physical constraints constitutes a novel framework. Hence, the focus of the current work is proof-of-concept and providing an investigation of the feasibility and capabilities of the proposed framework. The hitherto unknown properties of such physics-informed games, in combination with
the possibly very wide range of different outcomes depending on the physical problem (storage site properties), motivates an exploratory investigation. Therefore, we have performed a general mapping of the decision space by means of Pareto fronts and the use of a-posteriori decision-making methods. With maturing technology and if a clear decision-making procedure has already been established a-priori, the computational cost can be reduced by ignoring cases of limited interest. If that is indeed the case, there are multiple possibilities. For the special case of non-weighted utopia point Pareto front selection and social welfare maximization, the problem is significantly simplified and can be reduced to a single single-optimization problem. In the general case, however, multi-objective optimization needs to be performed. It is likely that by replacing the generic numerical methods used in this work (CSO and NSGA-II) by problem-adapted algorithms, significant computational cost reduction can be achieved.

In cooperative games, binding agreements are typically assumed only to be made between agents within the same coalitions. The wide extent of the Pareto fronts that have been observed in the numerical results here, implies that a wide range of coalition performances is possible. With the  potential disastrous effects on injection limits imposed by other coalitions and that correspond to the negative ends of the Pareto fronts, it is not likely that any coalition would just hope for the best when it comes to the actions of the external coalitions. Instead, it is far more sensible that distinct coalitions seek to reach agreements pertaining to all coalitions regarding injection/extraction limits before entering activity on a new site. Note that this does not mean that the grand coalition will necessarily form; coalition formation will still only occur when it is beneficial for agents to pool their resources and act as a single agent in terms of site activity and possibly transfer utility within the coalition. Hence, the range of scenarios indicated by the Pareto fronts are not meant to be interpreted as equally likely outcomes, but as an illustration of what could happen unless agreements and other regulations are in place.

Throughout the numerical test cases, we have equated the agent's payoff with its total injection/extraction contribution, i.e., the flexibility of transferable utility -- a very reasonable assumption for these games -- has not been investigated to its fullest extent. In reality, a coalition as a whole may benefit from some members giving up some of their value to another agent within the same coalition, i.e., allowing payoffs that are not strictly proportional to what a single agent injects/extracts. In this way, otherwise unstable coalition structures could be stabilized, for over all increased benefits.  

The objective function for any individual agent is arbitrary, so different transportation costs for different wells, relevant in CO$_2$ storage, can be directly incorporated into the objective functions. However, in the situation where the objectives are the values of injected fluid, and that value is not only dependent on the actual well and time, but also on the other agents of the same coalition structure, Theorem~\ref{thm:hierarch_PS} no longer applies. The proposed framework can still be employed, but the Pareto sets may not be subsets of each other. Hence, there may be no computational cost reduction associated with the hierarchical decomposition. 

Pareto fronts depicting the range of optimal injection strategies have been shown for test cases with specific setups. As demonstrated by the two  Bjarmeland formation test cases, different well configurations on the same site leads to different Pareto fronts. Different configurations of wells on a different subsurface site most likely yields very different Pareto fronts. Hence, one should be careful and avoid generalization of optimal injection schedules to other sites. Decision-making should always be based on site-specific models.

In this work, the evaluation of the constraints is more costly than evaluating the objective function. The number of constraint evaluations is equal to the number of objective function evaluation in the numerical methods, and the situation where the objective function is the more costly can be handled within the proposed numerical framework without any major changes.

A general concern in computational game theory is how to efficiently deal with the coalition structure generation problem with a large number of agents, especially in the presence of externalities of different sign that typically require searching the full space of possible coalition structures. In the current setting where the values of the coalition structures can only be obtained from numerical optimization and repeated evaluation of physical models, the computational cost quickly becomes infeasible with an increasing number of agents unless Theorem~\ref{thm:hierarch_PS} is applicable. Fortunately, even for the cases where the theorem does not apply, the number of agents on a realistic subsurface site is not likely to be more than a handful. The test cases with three and four agents we have presented here are considered representative with respect to the number of agents. The number of coalition structures with, e.g., five agents, is already quite large, but with increasing number of agents there may be domain-specific circumstances that makes it unnecessary to consider all possible collaborations, for instance between agents that for some reason are extremely unlikely to actually form a collaboration. As a consequence of Theorem~\ref{thm:hierarch_PS}, we can in principle obtain all coalition structure values from a sufficiently resolved Pareto front for the coalition structure with singleton coalitions only. With this approach, it is sufficient to solve a single multi-objective optimization problem, no matter the number of agents. Hence, the proposed framework is not restricted to academic test cases only.

\section{Conclusions}
We have shown that competition for subsurface resources can be modeled as a class of cooperative games in partition function form, where the coalition value (partition) function is a Pareto efficient solution to a multi-objective optimization problem. We have demonstrated that members of this class of games can display zero, negative, or both positive or negative externalities.  For certain value functions, there exist games with a unique solution identical for all coalition structures, inducing no preference on collaboration. If obtaining only this solution is the main goal, the computational problem is significantly simplified. If other solutions are of interest, the unique identical solution can still be used for reduced computational cost to obtain the full Pareto front.  

Furthermore, for a general class of games, we prove that the Pareto sets of the corresponding multi-objective optimization problems are hierarchically related so that the Pareto set of a coarser coalition structure is a proper subset of any coalition structure obtained by splitting coalitions of the coarser coalition structure. As demonstrated by numerical experiments, exploiting the hierarchical structure in the numerical methods to obtain the Pareto optimal solution sets leads to significantly reduced computational cost compared to a brute-force approach. The hierarchical property implies that all coalition structure values can be generated from the Pareto set corresponding to the singleton coalition structure. This indicates that problems with multiple agents may be solved efficiently 
since it is sufficient to compute a single Pareto front instead of a Pareto front for every possible coalition structure, which quickly becomes infeasible. The numerical results further display a wide range of  different games in terms of the values that can be obtained by the different coalitions. Hence, whether and how agents should collaborate is highly dependent on the specific game among the class of games introduced in this work, and conclusions should be game-specific before applied in practice.

\section*{Acknowledgments}
 This work was partly financed by the Norwegian Research Council grant 336294,
Expansion of Resources for CO$_2$ Storage on the Horda Platform (ExpReCCS). The authors would like to thank Auref Rostamian for suggestions regarding combined single- and multi-objective optimization.

\bibliographystyle{elsarticle-num}
\bibliography{bibliography} 

\end{document}